\documentclass[11pt,letterpaper]{compositio}

\usepackage[utf8]{inputenc}
\usepackage{amsmath}
\usepackage{amsthm}
\usepackage{amssymb}
\usepackage{tikz-cd}
\usepackage{cleveref}
\usepackage[shortlabels]{enumitem}
\usepackage{dynkin-diagrams}
\usepackage[notocbasic]{nomencl}
\makenomenclature

\newtheorem{thm}{Theorem}[section]
\newtheorem{theorem}[thm]{Theorem}
\newtheorem{corollary}[thm]{Corollary}

\newtheorem{lemma}[thm]{Lemma}
\newtheorem{prop}[thm]{Proposition}
\newtheorem{proposition}[thm]{Proposition}

\newtheorem{conjecture}[thm]{Conjecture}

\theoremstyle{definition}
\newtheorem{defn}[thm]{Definition}

\newtheorem{remark}[thm]{Remark}

\newcommand{\vfp}{V^{\mathrm{fp}}}
\newcommand{\vfploc}{V^{\mathrm{fp}}_{p(v)}}

\newcommand{\Aw}{\mathcal{A}_{w, \mathbb{F}_q}}
\newcommand{\AFq}{\mathcal{A}_{\mathbb{F}_q}}
\newcommand{\AwL}{\mathcal{A}_{w, \mathbb{F}_q}^{\mathcal{L}}}

\title[Polishchuk's conjecture and Kazhdan-Laumon representations]{Polishchuk's conjecture \\and \\Kazhdan-Laumon representations}
\author{Calder Morton-Ferguson}
\email{caldermf@stanford.edu}
\address{Mathematics Department, Stanford University\\450 Jane Stanford Way\\Stanford, CA 94305}
\date{\today}


\begin{document}

\begin{abstract}
    In their 1988 paper ``Gluing of perverse sheaves and discrete series representations," D. Kazhdan and G. Laumon constructed an abelian category $\mathcal{A}$ associated to a reductive group $G$ over a finite field with the aim of using it to construct discrete series representations of the finite Chevalley group $G(\mathbb{F}_q)$. The well-definedness of their construction depended on their conjecture that this category has finite cohomological dimension. This was disproven by R. Bezrukavnikov and A. Polishchuk in 2001, who found a counterexample in the case $G = SL_3$. Polishchuk then made an alternative conjecture: though this counterexample shows that the Grothendieck group $K_0(\mathcal{A})$ is not spanned by objects of finite projective dimension, he noted that a graded version of $K_0(\mathcal{A})$ can be thought of as a module over Laurent polynomials and conjectured that a certain localization of this module is generated by objects of finite projective dimension, and suggested that this conjecture could lead toward a proof that Kazhdan and Laumon's construction is well-defined. He proved this conjecture in Types $A_1, A_2, A_3$, and $B_2$. In the present paper, we prove Polishchuk's conjecture for all types, and prove that Kazhdan and Laumon's construction is indeed well-defined, giving a new geometric construction of discrete series representations of $G(\mathbb{F}_q)$.
\end{abstract}

\maketitle
\tableofcontents

\section{Introduction}

In their 1988 paper \cite{KL}, Kazhdan and Laumon described a gluing construction for perverse sheaves on the basic affine space associated to a semisimple algebraic group $G$ split over a finite field $\mathbb{F}_q$, defining an abelian category $\mathcal{A}$ of ``glued perverse sheaves" consisting of certain tuples of perverse sheaves on the basic affine space indexed by the Weyl group. They aimed to use these categories to provide a new geometric construction of discrete series representations of $G(\mathbb{F}_q)$.

Their proposal was to use $\mathcal{A}$ to construct representations as follows. First, they observed that the discrete series representations they sought to construct arise from characters of the non-split tori $T(w)$ of $G$, which are indexed by the Weyl group. For each $w \in W$, they defined a category $\mathcal{A}_{w, \mathbb{F}_q}$ in a way such that $K_0(\mathcal{A}_{w, \mathbb{F}_q})$ carries commuting actions of $G(\mathbb{F}_q)$ and $T(w)$.

They expected that the wildly infinite-dimensional representation $$K_0(\mathcal{A}_{w, \mathbb{F}_q}) \otimes \mathbb{C}$$ of $G(\mathbb{F}_q)$ admits a finite-dimensional quotient whose $T(w)$-isotypic components are the discrete series representations they sought to construct. Following the philosophy of Grothendieck's sheaf-function dictionary, Kazhdan and Laumon knew that the appropriate subspace of $K_0(\mathcal{A}_{w, \mathbb{F}_q}) \otimes \mathbb{C}$ by which one should take the quotient should be the kernel of a certain ``Grothendieck-Lefschetz pairing" on $K_0(\mathcal{A}_{w, \mathbb{F}_q})$, which is defined in terms of the $\mathrm{Ext}$ groups in the category $\mathcal{A}$. They then made the following conjecture and proved that it implies the well-definedness of their representations.

\begin{conjecture}[(Kazhdan-Laumon, \cite{KL})]\label{conj:kl}
    The category $\mathcal{A}$ has finite cohomological dimension. In other words, for any two objects $A$ and $B$ of $\mathcal{A}$, there is an $n$ for which $\mathrm{Ext}^i(A, B) = 0$ whenever $i > n$.
\end{conjecture}
More than a decade later, Bezrukavnikov and Polishchuk found a counterexample to this conjecture in the case $G = \mathrm{SL}_3$.

\begin{proposition}[(Bezrukavnikov-Polishchuk, Appendix to \cite{P})]\label{prop:bp}
    \Cref{conj:kl} is false.
\end{proposition}

In \cite{P}, Polishchuk put forward the idea that although \Cref{conj:kl} is false as stated in \cite{KL}, it is not strictly necessary in order to prove the more important assertion that Kazhdan and Laumon's construction of representations is well-defined. He notes that $K_0(\mathcal{A}_{w, \mathbb{F}_q})$ carries the structure of a $\mathbb{Z}[v, v^{-1}]$-module using the formalism of mixed sheaves where $v$ acts by a Tate twist, and then frames \Cref{conj:kl} as the claim that $K_0(\mathcal{A}_{w, \mathbb{F}_q})$ is spanned by objects of finite projective dimension. In this situation, the Grothendieck-Lefschetz pairing defined on $K_0(\mathcal{A}_{w, \mathbb{F}_q}) \otimes \mathbb{C}$ can be thought of as taking polynomial values in $\mathbb{Z}[v, v^{-1}]$ and then specializing at $v = q^{\frac{1}{2}}$, which is one way to see why \Cref{conj:kl} would imply the well-definedness of this pairing and therefore the well-definedness of Kazhdan and Laumon's construction.

Although \Cref{prop:bp} shows that this is false, he instead proposes that this pairing is still well-defined if one allows it to take values in a certain localization of the ring $\mathbb{Z}[v, v^{-1}]$. Letting $\mathcal{A}_{\mathbb{F}_q} = \mathcal{A}_{e, \mathbb{F}_q}$, i.e. the category of ``Weil sheaves'' in the Kazhdan-Laumon context, Polishchuk proposes the following more precise conjecture as a first step toward this goal.

\begin{conjecture}\label{conj:p}
    There exists a finite set of polynomials, which are nonzero away from roots of unity, such that the localization of $K_0(\AFq)$ at the multiplicative set generated by these polynomials is generated by objects of finite projective dimension.
\end{conjecture}

In \cite{P}, Polishchuk develops a framework toward answering this conjecture, resolving it himself in Types $A_1, A_2, A_3,$ and $B_2$. In our first main theorem, we use this framework along with the algebraic understanding of symplectic Fourier transforms provided by \cite{CMFSymplectic} to prove this conjecture in general.
\begin{theorem}\label{thm:polyconj}
    \Cref{conj:p} is true. In particular, the localization of the $\mathbb{Z}[v, v^{-1}]$-module $K_0(\AFq)$ at the polynomial
    \begin{align}
        p(v) & = \prod_{i=1}^{\ell(w_0)} \left(1 - v^{2i}\right)
    \end{align}
    is generated by objects of finite projective dimension.
\end{theorem}

As Polishchuk expected, the resolution of \Cref{conj:p} brings us very close to showing that Kazhdan and Laumon's construction of discrete series representations is well-defined. The main result of the present paper is that by using the formalism of monodromic perverse sheaves, we can prove a similar theorem which indeed completes the necessary technicalities to carry out Kazhdan and Laumon's construction in general.
\begin{theorem}\label{thm:mainthm}
    For any character sheaf $\mathcal{L}$ of $T$ and element $w \in W$, the localization of the $\mathbb{Z}[v, v^{-1}]$-module $K_0(\AwL)$ at $p(v)$ is spanned by classes of objects of finite projective dimension in $\AwL$.
\end{theorem}

This monodromic approach to Kazhdan and Laumon's construction was already successfully carried out in \cite{BP}, in which Braverman and Polishchuk explain how to carry out a well-defined version of Kazhdan and Laumon's construction in the case where $\mathcal{L}$ corresponds to a \emph{quasi-regular} character. So one can think of the following corollary to Theorem \ref{thm:mainthm} as a generalization of Braverman and Polishchuk's result to the case of an arbitrary character.
\begin{corollary}
    The Kazhdan-Laumon construction proposed in \cite{KL} is well-defined for monodromic sheaves corresponding to any character.
\end{corollary}

\subsection{Layout of the paper} In Section \ref{sec:prelim}, we explain some background on Kazhdan-Laumon categories. In Section \ref{sec:monodrom}, we then explain the monodromic setting required to state \Cref{thm:mainthm} and discuss the categorical center of the monodromic Hecke category, which will be an important tool in the proof. Then in Section \ref{sec:dg} we use dg formalism to explain why the derived category of the Kazhdan-Laumon category admits an action of this categorical center. This is followed in Section \ref{sec:canonical} by an explanation of a crucial tool in the study of Kazhdan-Laumon categories proposed by Polishchuk \cite{P} called the canonical complex. We then complete the proofs of our results: in Section \ref{sec:central} we prove \Cref{thm:mainthm}, and in Section \ref{sec:final} we recall the results of \cite{CMFSymplectic} and explain how it, combined with the previous setup, allow us to prove Polishchuk's original conjecture and establish \Cref{thm:polyconj} independently of the monodromic setting. Finally, in \Cref{sec:construction}, we explain how to carry out the construction of Kazhdan-Laumon representations explicitly given our theorems.

\subsection{Acknowledgments}
I am very grateful for the support of my advisor, Roman Bezrukavnikov, who introduced me to Kazhdan-Laumon's construction and provided continuous feedback and support throughout all stages of the project. I would also like to thank Alexander Polishchuk, Pavel Etingof, Zhiwei Yun, Ben Elias, Minh-Tam Trinh, Elijah Bodish, Alex Karapetyan and Matthew Nicoletti for helpful conversations. During this work, I was supported by
an NSERC PGS-D award.

\section{Preliminaries}\label{sec:prelim}

\subsection{Background and notation}

\subsubsection{General setup} Let $G$ be a split semisimple group over a finite field $\mathbb{F}_q$. Let $T$ be a Cartan subgroup split over $\mathbb{F}_q$, $B$ a Borel subgroup containing $T$, and $U$ its unipotent radical. Let $X = G/U$ be the basic affine space associated to $G$ considered as a variety over $\mathbb{F}_q$. Let $W$ be the Weyl group $W$. We let $S$ denote the set of simple reflections in $W$. Writing $q = p^m$ for some prime number $p$, we choose $\ell$ to be a prime with $\ell \neq p$.

\subsubsection{$\ell$-adic sheaves, Tate twists, and Grothendieck groups}\label{sec:ladic}

We work with the category of perverse sheaves $\mathrm{Perv}(G/U)$ of mixed $\ell$-adic perverse sheaves on the basic affine space $G/U$, and more generally with the constructible derived category $D^b(G/U)$ of mixed $\ell$-adic sheaves on $G/U$. We choose an isomorphism $\overline{\mathbb{Q}}_{\ell} \cong \mathbb{C}$ and work with $\mathbb{C}$ going forward. Pick a square root $q^{\frac{1}{2}}$ of $q$ in $\mathbb{C}$ once and for all, and define the half-integer Tate twist $(\tfrac{1}{2})$ on $D^b(G/U)$. 
We then view $K_0(G/U) = K_0(D^b(G/U))$ as a $\mathbb{Z}[v, v^{-1}]$-module where $v^{-1}$ acts by $(\tfrac{1}{2})$. When $\mathcal{C}$ is any category for which $K_0(\mathcal{C})$ is a $\mathbb{Z}[v, v^{-1}]$-module, we denote by $K_0(\mathcal{C})\otimes \mathbb{C}$ the specialization $K_0(\mathcal{C}) \otimes_{\mathbb{Z}[v, v^{-1}]} \mathbb{C}$ at $v = q^{\frac{1}{2}}$. We use $\mathbb{D}$ to denote the Verdier duality functor. We let ${}^pH^i$ be the perverse cohomology functors for any $i \in \mathbb{Z}$. 

We also choose once and for all a nontrivial additive character $\psi : \mathbb{F}_q \to \overline{\mathbb{Q}_{\ell}}$, and let $\mathcal{L}_{\psi}$ be the corresponding Artin-Schrieier sheaf on $\mathbb{G}_a$.

The variety $G/U$ comes with a natural Frobenius morphism $\mathrm{Fr} : G/U \to G/U$; we can then consider the category of Weil sheaves $\mathrm{Perv}_{\mathbb{F}_q}(G/U)$ on $G/U$, i.e. sheaves $K \in \mathrm{Perv}(G/U)$ equipped with a natural isomorphism $\mathrm{Fr}^*K \cong K$.

\subsubsection{Elements of $G$ indexed by $W$} For every simple root $\alpha_s$ of $G$ corresponding to the simple reflection $s$, we fix an isomorphism of the corresponding one-parameter subgroup $U_s \subset U$ with the additive group $\mathbb{G}_{a}$. This uniquely defines a homomorphism $\rho_s : SL_{2} \to G$ which induces the given isomorphism of $\mathbb{G}_{a}$ (embedded in $SL_{2}$ as upper-triangular matrices) with $U_s$; then let
\begin{align*}
    n_s & = \rho_s
    \begin{pmatrix}
        0 & 1\\-1 & 0 
    \end{pmatrix}.
\end{align*}
For any $w \in W$, writing a reduced word $w = s_{i_1} \dots s_{i_k}$ we set $n_w = n_{s_{i_1}} \dots n_{s_{i_k}}$, and one can check that this does not depend on the reduced word. We also define for any $s \in S$ the subtorus $T_s \subset T$ obtained from the image of the coroot $\alpha_s^\vee$ and define $T_w$ for any $w \in W$ to be the product of all $T_s$ ($s \in S$) for which $s \leq w$ in the Bruhat order.

\subsection{Kazhdan-Laumon categories}

\subsubsection{Fourier transforms on $\mathrm{Perv}(G/U)$}\label{sec:ft}

In \cite{KL} and \cite{P}, to each $w \in W$ the authors assign an element of $D^b(G/U\times G/U)$ which, up to shift, is perverse and irreducible. Following \cite{P}, let $X(w) \subset G/U \times G/U$ be the subvariety of pairs $(gU, g'U) \subset (G/U)^2$ such that $g^{-1}g' \in Un_wT_wU$. There is a canonical projection $\mathrm{pr}_w : X(w) \to T_w$ sending $(gU, g'U)$ to the unique $t_w \in T_w$ such that $g^{-1}g' \in Un_wt_wU$. In the case when $w = s \in S$, the morphism $\mathrm{pr}_s : X(s) \to T_s \cong \mathbb{G}_{m, k}$ extends to $\overline{\mathrm{pr}}_s : \overline{X(s)} \to \mathbb{G}_{a,k}$ and we have
\begin{align*}
    \overline{K(s)} = (-\overline{\mathrm{pr}}_s)^* \mathcal{L}_\psi
\end{align*}
and in the case of general $w \in W$
\begin{align}\label{eqn:wfroms}
    \overline{K(w)} & = \overline{K(s_{i_1})} * \dots * \overline{K(s_{i_k})}
\end{align}
whenever $w = s_{i_1} \dots s_{i_k}$ is a reduced expression, where $*$ denotes the convolution of sheaves on $G/U \times G/U$ as defined in \cite{KL}. One can take this as the definition of Kazhdan-Laumon sheaves, which is well-defined by the proposition below, or refer to the explicit definition of $\overline{K(w)}$ which works for all $w \in W$ at once given in \cite{KL} or \cite{P}.

\begin{proposition}[\cite{KL}]
    The Kazhdan-Laumon sheaves $\overline{K(s)}$ for $s \in S$ under convolution satisfy the braid relations (up to isomorphism).
\end{proposition}

For any $s \in S$, they note that the endofunctor $K \to K * \overline{K(s)}$ of $D^b(G/U)$ can be identified with a certain ``symplectic Fourier-Deligne transform'' defined as follows. Let $P_s$ be the parabolic subgroup of $G$ associated to $s$, and let $Q_s = [P_s, P_s]$. The map $G/U \to G/Q_s$ has all fibers isomorphic to $\mathbb{A}^2 \setminus \{(0, 0)\}$, and it is shown in Section 2 of \cite{KL} that there exists a natural fiber bundle $\pi : V_s \to G/Q_s$ of rank $2$ equipped with a $G$-invariant symplectic pairing which contains $G/U$ as an open set, with inclusion $j : G/U \to V_s$ with $\pi\circ j$ being the original projection $G/U \to G/Q_s$. There is then a symplectic Fourier-Deligne transform $\tilde{\Phi}_s$ on $D^b(V_s)$ defined by
\begin{align}
    \tilde{\Phi}_s(K) & = p_{2!}(\mathcal{L}\otimes p_1^*(K))[2](1),
\end{align}
where the $p_i$ are the projections of the product $V_s \times_{G/Q_s} V_s$ on its factors, and $\mathcal{L} = \mathcal{L}_\psi(\langle, \rangle)$ is a smooth rank-1
$\overline{\mathbb{Q}_\ell}$-sheaf which is the pullback of the Artin—Schreier sheaf $\mathcal{L}_\psi$ under the morphism $\langle, \rangle$, c.f. Section 4 of \cite{P}. We then define the endofunctor $\Phi_s$ of $D^b(G/U)$ By
\begin{align}
    \Phi_s(K) & = j^*\tilde{\Phi}_sj_{!} K.
\end{align}

\begin{prop}[\cite{KL}, \cite{P}]
    The functors $\Phi_s$ and $- * \overline{K(s)}$ are naturally isomorphic.
\end{prop}

For any $w \in W$, we let $\Phi_w = \Phi_{s_{i_1}} \dots \Phi_{s_{i_k}}$ where $s_{i_1} \dots s_{i_k}$ is a reduced expression for $w$ as a product of simple reflections. The functors $\Phi_w$ are the gluing functors which Kazhdan and Laumon use to define the so-called glued categories $\mathcal{A}$.

\begin{defn}
    Using the six-functor formalism, one can check that each $\Phi_s$ has a right adjoint which we call $\Psi_s$, following the setup of Section 1.2 of \cite{P}. The functors $\Psi_s$ also form a braid action on $D^b(G/U)$. We then define $\Psi_w$ similarly.
\end{defn}
The functors $\Phi_w$ (resp. $\Psi_w$) are each right (resp. left) $t$-exact on $D^b(G/U)$ with respect to the perverse $t$-structure. For any $w \in W$, let $\Phi_w^\circ = {}^pH^0\Phi_w$ and $\Psi_w^\circ = {}^pH^0\Psi_w$, noting that $\Phi_w = L\Phi_w^\circ$, $\Psi_w = R\Psi_w^\circ$.

\begin{prop}[\cite{P}, Section 4.1]
    For any $s \in S$, there are natural morphisms $c_s : \Phi_s^2 \to \mathrm{Id}$ and $c_s' : \mathrm{Id} \to \Psi_s^2$ satisfying the associativity conditions
    \begin{align}
        \Phi_s \circ c_s & = c_s \circ \Phi_s : \Phi_s^3 \to \Phi_s & \Psi \circ c_s' & = c_s' \circ \Psi : \Psi_s \to \Psi_s^3.
    \end{align}
\end{prop}

\begin{corollary}\label{cor:nu}
    For any $y', y \in W$, there is a natural transformation $\nu_{y', y} : \Phi_{y'}\Phi_y \to \Phi_{y'y}$.
\end{corollary}
\begin{proof}
    We go by induction on $\ell(y) + \ell(y')$. If $\ell(y'y) = \ell(y') + \ell(y)$, then $\nu_{y', y}$ is the tautological map arising from the fact that the $\Phi_w$ form a braid action. If instead $\ell(y'y) < \ell(y') + \ell(y)$, then there exists some $s \in S$ such that $y' = \tilde{y}'s$ and $y = s\tilde{y}$ for some $\tilde{y}', \tilde{y} \in W$ with $\ell(\tilde{y}'s) = \ell(\tilde{y}') + 1$ and $\ell(s\tilde{y}) = \ell(\tilde{y}) + 1$, and so we have maps
    \[\begin{tikzcd}
        \Phi_{y'}\Phi_y = \Phi_{\tilde{y}'}\Phi_s^2\Phi_{\tilde{y}} \arrow[rr, "\Phi_{\tilde{y}}\circ c_s \circ \Phi_{\tilde{y}}"] & & \Phi_{\tilde{y}'}\Phi_{\tilde{y}} \arrow[r, "\nu_{\tilde{y}',\tilde{y}}"] & \Phi_{y'y},
    \end{tikzcd}\]
    the former coming from $c_s$ and the latter coming from our induction hypothesis.
\end{proof}

\subsubsection{Definition of the Kazhdan-Laumon category}

\begin{defn}[\cite{KL}, \cite{P}]\label{def:kldef}
    The Kazhdan-Laumon category $\mathcal{A}$ has objects which are tuples $(A_w)_{w \in W}$ with $A_w \in \mathrm{Perv}(G/U)$ and equipped with morphisms 
    \begin{align}
        \theta_{y,w} : \Phi_y^\circ A_w \to A_{yw}
    \end{align}
    for every $y, w \in W$ such that the diagram
    \[\begin{tikzcd}
        \Phi_{y'}^\circ\Phi_{y}^\circ A_w \arrow[r, "\Phi_{y'}\theta_{y,w}"] \arrow[d, "\nu_{y',y}"] & \Phi_{y'}^\circ A_{yw} \arrow[d, "\theta_{y',yw}"]\\
        \Phi_{y'y}^\circ A_w \arrow[r, "\theta_{y'y,w}"] & A_{y'yw}
    \end{tikzcd}\]
    commutes for any $y, y', w \in W$. 

    A morphism $f$ between objects $(A_w)_{w \in W}$ and $(B_w)_{w \in W}$ is a collection of morphisms $f_w : A_w \to B_w$ such that
    \[\begin{tikzcd}
        \Phi_y^\circ A_w \arrow[r, "\Phi_y^\circ f_w"] \arrow[d, "\theta_{y,w}^A"] & \Phi_y^\circ B_w\arrow[d, "\theta_{y,w}^B"]\\
        A_{yw} \arrow[r, "f_{yw}"] & B_{yw}.
    \end{tikzcd}\]
\end{defn}
It is shown in \cite{P} that this category is abelian, and that the functors $j_{w}^* : \mathcal{A} \to \mathrm{Perv}(G/U)$ defined by $j_{w}^*((A_w)_{w \in W}) = A_w$ are exact.

\begin{remark}
    We could have instead asked for morphisms $A_{yw} \to \Psi_y^\circ A_w$, making reference to the functors $\Psi_y^\circ$ rather than the $\Phi_y^\circ$. Later, we will discuss an alternate and more elegant definition of $\mathcal{A}$ as coalgebras over a certain comonad on $\oplus_{w \in W} \mathrm{Perv}(G/U)$ which is assembled from the functors $\Psi_y^\circ$. In \Cref{def:kldef} though, we present the definition of the Kazhdan-Laumon category as it was originally formulated in \cite{KL} and later explained in more detail in \cite{P}.
\end{remark}

\subsubsection{Definition of the $w$-twisted categories $\mathcal{A}_{w, \mathbb{F}_q}$}

We note that the category $\mathcal{A}$ carries an action of the Weyl group $W$ as follows. For any $w \in W$, we let $\mathcal{F}_w : \mathcal{A} \to \mathcal{A}$ be the exact functor defined by right translation of the indices in the tuple, i.e. $\mathcal{F}_w((A_y)_{y\in W}) = (A_{yw})_{y \in W}$.

\begin{defn}\label{def:aw}
    For any $w \in W$, let $\mathcal{A}_{w,\mathbb{F}_q}$ be the category with objects $(A, \psi_A)$ where $A \in \mathcal{A}$ and $\psi_A : \mathcal{F}_w\mathrm{Fr}^*A \to A$ is an isomorphism. We call these $w$-twisted Weil sheaves in the Kazhdan-Laumon category.

    For any two such objects $(A, \psi_A)$ and $(B, \psi_B)$, we let $\mathrm{Hom}_{\mathcal{A}_{w, \mathbb{F}_q}}(A, B)$ be the set of morphisms $f \in \mathrm{Hom}_{\mathcal{A}}(A, B)$ such that $f \circ \psi_A = \psi_B \circ \mathcal{F}_w\mathrm{Fr}^*f$.
\end{defn}

\begin{remark}
    When $w = e$ is trivial, we will write $\mathcal{A}_{\mathbb{F}_q} = \mathcal{A}_{w, \mathbb{F}_q}$. In this case, it is shown in \cite{P} that $\mathcal{A}_{\mathbb{F}_q}$ is equivalent to the category obtained by applying the Kazhdan-Laumon gluing procedure described in \Cref{def:kldef} to the category $\mathrm{Perv}_{\mathbb{F}_q}(G/U)$ of Weil perverse sheaves on $G/U$, i.e. perverse sheaves $K$ equipped with an isomorphism $\mathrm{Fr}^*K \to K$.
\end{remark}

\subsubsection{Adjoint functors on $D^b(\mathcal{A})$}
\begin{defn}
    For any $w \in W$, let $j_w^* : D^b(\mathcal{A}) \to D^b(G/U)$ be the functor arising from the same-named exact functor $j_w^* : \mathcal{A} \to \mathrm{Perv}(G/U)$ given by $j_w^*((A_y)_{y \in W}) = A_w$.

    We then define a functor $j_{w!}^\circ : \mathrm{Perv}(G/U) \to \mathcal{A}$ by
    \begin{align}
        j_{w!}^\circ(K) & = (\Phi_{yw^{-1}}^\circ K)_{y \in W}.
    \end{align}
    One can check that the morphisms $\nu_{y', y}$ for $y, y' \in W$ introduced in \Cref{cor:nu} endow the tuple $(\Phi_{yw^{-1}}^\circ K)_{y \in W}$ with the structure morphisms required to define an object of $\mathcal{A}$. We let $j_{w!}$ be the left-derived functor to $j_{w!}^\circ$.
\end{defn}

\begin{prop}[(Proposition 7.1.2, \cite{P})]\label{prop:adjunctionshriek}
    For any $w \in W$, there is an adjunction $(j_{w!}^\circ, j_w^*)$. Further, the functor $j_{w!} : D^b(G/U) \to D^b(\mathcal{A})$ has the property that
    \begin{align}
        {}^pH^i(j_{w!}(K)) & = ({}^pH^i\Phi_{yw^{-1}}(K))_{y \in W},\label{eqn:jwshriek}
    \end{align}
    and there is also an adjunction $(j_{w!}, j_w^*)$ on derived categories.
\end{prop}
Analogously, acting instead by the functors $\Psi_{yw^{-1}}$ in (\ref{eqn:jwshriek}) defines a right-adjoint $j_{w*}^\circ$ to $j_w^*$ and its right-derived functor $j_{w*}$ in the very same way, as is also explained in \cite{P}.

\subsubsection{The functor $\iota$}

\begin{defn}
    We define an endofunctor $\iota$ of $\mathcal{A}$ by 
    \begin{align}
        \iota((A_w)_{w \in W}) & = (\Phi_{w_0}^\circ A_{w_0w})_{w \in W}
    \end{align}
    with the structure morphisms described in 7.2 of \cite{P}. 

    It is shown in loc. cit. that we can also abuse notation and view $\iota$ as a functor on $D^b(\mathcal{A})$ (by replacing $\Phi_{w_0}^\circ$ with $\Phi_{w_0}$); for our purposes we will only need the functor $\iota^2$, which we will later describe as an endofunctor of $D^b(\mathcal{A})$ more conceptually in \Cref{lem:i2ft}.
\end{defn}

\subsubsection{Objects of the form $j_{w!}(A)$}

\begin{prop}\label{prop:adjunction}
    For any $B \in D^b(G/U)$ and any $w \in W$, the object $j_{w!}(B) \in D^b(\mathcal{A})$ has finite projective dimension. 
\end{prop}

\begin{proof}
    The adjunction in Proposition \ref{prop:adjunctionshriek} gives that for any $A \in \mathcal{A}$,
    \begin{align}
        \mathrm{Ext}_{\mathcal{A}}^\bullet(j_{w!}(B), A) & = \mathrm{Ext}_{\mathrm{Perv}(G/U)}^\bullet(B, j_w^*A),
    \end{align}
    and the latter is a finite-dimensional vector space because the category $\mathrm{Perv}(G/U)$ has finite cohomological dimension.
\end{proof}

\subsubsection{Polynomials and localization}
Fix the Weyl group $W$ and choose $w_0$ its longest element.

\begin{defn}
    Define $P(x, v), \tilde{P}(x, v) \in \mathbb{Z}[x, v, v^{-1}]$ by
    \begin{align}
        P(x, v) & = \prod_{i=0}^{\ell(w_0)} (x - v^{2i})\\
        \tilde{P}(x, v) & = \prod_{i=1}^{\ell(w_0)} (x - v^{2i})
    \end{align}
    and let $p(v) = \tilde{P}(1, v)$.
\end{defn}

\begin{defn}
    For any abelian category $\mathcal{C}$ such that $K_0(\mathcal{C})$ has the structure of a $\mathbb{Z}[v, v^{-1}]$-module, let $\vfp \subset K_0(\mathcal{C})$ be the submodule spanned by all objects of finite projective dimension. Further, let $\vfploc$ denote the localization of this module at $p(v)$, or equivalently at all of the linear factors $(1 - v^{2i})$ for $1 \leq i \leq \ell(w_0)$.
\end{defn}

\section{Monodromic sheaves, character sheaves, and the categorical center}\label{sec:monodrom}

\subsubsection{Rank one character sheaves on $T$}

We let $\mathrm{Ch}(T)$ be the category of rank one character sheaves on $T$; we refer to Appendix A of \cite{YCh} for a detailed treatment. We note that $\mathrm{Ch}(T)$ carries a natural action of the Weyl group $W$.

For any $\mathcal{L}$ in $\mathrm{Ch}(T)$, we let $W_{\mathcal{L}}^\circ$ be the normal subgroup of the stabilizer of $\mathcal{L}$ in $W$ which is the Weyl group of the root subsystem of the root system of $W$ on which $\mathcal{L}$ is trivial; see 2.4 of \cite{lusztig_yun_2020} for details.

\subsubsection{Monodromic version of the Kazhdan-Laumon category}

In Sections 2.1.3 and 2.3.4 of \cite{BP}, the authors explain how to define a category $\mathrm{Perv}_{\mathcal{L}}(G/U)$ of monodromic sheaves on $G/U$ with respect to the monodromy $\mathcal{L}$. We then let $D_{\mathcal{L}}^b(G/U) = D^b(\mathrm{Perv}_{\mathcal{L}}(G/U))$ be its derived category.

\begin{remark}
    We note that in this definition, monodromic perverse sheaves are defined such that the extension of two $\mathcal{L}$-equivariant sheaves may not be $\mathcal{L}$-equivariant, only $\mathcal{L}$-monodromic. In other words, for $\mathcal{L}$ trivial, this reduces to the category of perverse sheaves on $G/U$ with unipotent monodromy (c.f. \cite{BY}) on the right rather than simply to $\mathrm{Perv}(G/B)$.
\end{remark}
We use this to define the $\mathcal{L}$-monodromic Kazhdan-Laumon category. First, it is straightforward to check that if $A \in \mathrm{Perv}_{\mathcal{L}}(G/U)$, then for any $w \in W$, $\Phi_w^\circ A \in \mathrm{Perv}_{w\mathcal{L}}(G/U)$ (c.f.\ \Cref{prop:ftconv}). We can then formulate the following definition.
\begin{defn}
    For any $\mathcal{L} \in \mathrm{Ch}(T)$, we define the \emph{monodromic Kazhdan-Laumon category} $\mathcal{A}^{\mathcal{L}}$ as the category whose objects are $(A_w)_{w \in W}$ with $A_w \in \mathrm{Perv}_{w\mathcal{L}}(G/U)$, and equipped with the same morphisms and compatibilities which appear in \Cref{def:kldef}. 
    
    If $\mathcal{L}$ is such that $\mathrm{Fr}^*(w\mathcal{L}) \cong \mathcal{L}$, we can then further define $\AwL$ in analogy to \Cref{def:aw} as the category of pairs $(A, \psi_A)$ where $A \in \mathcal{A}^{\mathcal{L}}$ and $\psi_A : \mathcal{F}_w\mathrm{Fr}^*A \to A$ an isomorphism.
\end{defn}

\subsection{The monodromic Hecke category and its center}

\subsubsection{The monodromic Hecke category}\label{sec:monodromhecke} In \cite{Gouttard}, the author defines for any $\mathcal{L}$ a category $\mathcal{P}_{\mathcal{L}}$ (called $D^b_{(B)}(G/U)_{t}$ in loc. cit., where $t$ is a parameter determined by $\mathcal{L}$) such that for $\mathcal{L}$ trivial, this reduces to the familiar derived category $D^b_{(B)}(G/U)$ of $B$-constructible sheaves on $G/U$. Elements of $\mathcal{P}_{\mathcal{L}}$ are $\mathcal{L}$-monodromic with respect to the right action of $T$, while their left monodromies may correspond to any character sheaf in the $W$-orbit of $\mathcal{L}$.

For any $w \in W$, this category contains monodromic versions ${}_{w\mathcal{L}}\Delta(w)_{\mathcal{L}}$ and ${}_{w\mathcal{L}}\nabla(w)_{\mathcal{L}}$ of standard and costandard sheaves; we emphasize that extensions of such objects in this category may not be $\mathcal{L}$-equivariant with respect to the right $T$-action even when they remain $\mathcal{L}$-monodromic, so this category also contains monodromic versions of tilting objects $\mathcal{T}(w)_{\mathcal{L}}$ for any $w \in W$.

We let ${}_{\mathcal{L}}\mathcal{P}_{\mathcal{L}}$ be the triangulated subcategory of objects generated by the standard and costandard objects corresponding to $w \in W_{\mathcal{L}}^\circ$, each of whose objects must also be left-monodromic with respect to $\mathcal{L}$.

\subsubsection{Free-monodromic Hecke categories}\label{subsubsec:free-monodromic}

In \cite{BY}, in the unipotent monodromy case where $\mathcal{L}$ is trivial, the authors define a category formed from a certain completion of $D^b_{\mathcal{L}}(G/U)$ called the category of \emph{unipotently free-monodromic} sheaves.

In \cite{Gouttard}, the case of non-unipotent monodromy was treated carefully. In loc. cit., the author defines a category $\hat{\mathcal{P}}_{\mathcal{L}}$ (which is called $\hat{D}^b_{(B)}(G/U)_{t}$ in loc. cit. where $t$ is a parameter determined by $\mathcal{L}$), which is a certain completion of the category $\mathcal{P}_{\mathcal{L}}$ defined in \Cref{sec:monodromhecke}, equipped with a monoidal structure which we also denote by $*$ in this context.

This category contains objects $\varepsilon_{n,\mathcal{L}}$ and $\hat{\delta}_{\mathcal{L}}$ introduced in Corollary 5.3.3 of \cite{BT}. The object $\hat{\delta}_{\mathcal{L}}$ is the monoidal unit for the convolution product on $\hat{\mathcal{P}}_{\mathcal{L}}$, whereas convolution with the objects $\varepsilon_{n, \mathcal{L}}$ can be thought of as a sort of projection to the subcategory of objects whose corresponding ``logarithmic monodromy operator" is nilpotent of order at most $n$; c.f. Appendix A of \cite{BY} for an explanation of this perspective.

Finally, we recall that $\hat{\mathcal{P}}_{\mathcal{L}}$ contains for all $w \in W$ free-monodromic versions ${}_{w\mathcal{L}}\hat{\Delta}(w)_{\mathcal{L}}$, ${}_{w\mathcal{L}}\hat{\nabla}(w)_{\mathcal{\mathcal{L}}}$ of standard and costandard sheaves, c.f. \cite{BY} for the unipotent case, \cite{lusztig_yun_2020} for a description of these objects in $\mathcal{P}_{\mathcal{L}}$ for arbitrary $\mathcal{L}$, and \cite{Gouttard} for their free-monodromic versions in $\hat{\mathcal{P}}_{\mathcal{L}}$. We define the subcategory ${}_{\mathcal{L}}\hat{\mathcal{P}}_{\mathcal{L}}$ analogously to the subcategory ${}_{\mathcal{L}}\mathcal{P}_{\mathcal{L}}$ of $\mathcal{P}_{\mathcal{L}}$.

\subsubsection{The center of the monodromic Hecke category}

To define the notion of categorical center which will be useful in this setting, we will closely follow the conventions of \cite{BITV} in this section. We begin by recalling some definitions from loc. cit.

\begin{defn}[\cite{BITV}]
    Let $\mathcal{Y} = (G/U \times G/U)/T$ where $T$ acts by the right-diagonal multiplication, and let $\mathcal{H}^{(1)} = D^b(G\backslash \mathcal{Y})$. Then for any $\mathcal{L} \in \mathrm{Ch}(T)$, let $\mathcal{H}^{(1)}_{\mathcal{L}}$ be the full $\mathcal{L}$-monodromic subcategory of $\mathcal{H}^{(1)}$ with respect to the projection $\mathcal{Y} \to (G/B)^2$.

    For $\mathfrak{o}$ a $W$-orbit in $\mathrm{Ch}(T)$, we let $\mathcal{H}^{(1)}_{\mathfrak{o}}$ be the full subcategory of $\mathcal{H}^{(1)}$ consisting of monodromic sheaves with monodromies in $\mathfrak{o}$, i.e.
    \begin{align}
        \mathcal{H}^{(1)}_{\mathfrak{o}} \cong \bigoplus_{\mathcal{L} \in \mathfrak{o}} \mathcal{H}^{(1)}_{\mathcal{L}}.
    \end{align}
\end{defn}
Following 3.2 of \cite{BITV}, consider the diagram
\[\begin{tikzcd}
    & G \arrow[dl, "\pi"'] \arrow[dr, "q"] \times G/B & \\
    G & & \mathcal{Y}
\end{tikzcd}\]
where $\pi$ is the projection and $q$ is the quotient of the map $q' : G \times G/U \to G/U \times G/U$ given by $q'(g, xU) = (xU, gxU)$ by the free right $T$-action, with respect to which $q'$ is equivariant. 
\begin{defn}
    The Harish-Chandra transform is the functor
    \begin{align}
        \mathfrak{hc} = q_! \circ \pi^*: D^b(G/_{\mathrm{Ad}} G) \to D^b(G\backslash \mathcal{Y}),
        \end{align}
        which is monoidal with respect to the natural convolution product on each side by \cite{GinsburgAdmissible}, c.f. 2.2 of \cite{BITV} for a detailed exposition of these convolution products.
\end{defn}

\begin{defn}\label{def:charsheaves}
    For any $\mathcal{L} \in \mathrm{Ch}(T)$, let $D^b_{\mathfrak{C},\mathcal{L}}(G) \subset D^b(G/_{\mathrm{Ad}}G)$ be the full triangulated subcategory with objects $\mathcal{F}$ satisfying $\mathfrak{hc}(\mathcal{F}) \in \mathcal{H}^{(1)}_{\mathfrak{o}}$, where $\mathfrak{o}$ is the $W$-orbit of $\mathcal{L}$.
\end{defn}

In \cite[Theorem 5.3.4]{BITV}, the authors show that the functor $\mathfrak{hc}$ realizes the category $D^b_{\mathfrak{C}, \mathcal{L}}(G)$ as the center $\mathcal{Z}\mathcal{H}^{(1)}_{\mathfrak{o}}$, with the notion of center in this case being defined in \cite[Section 2.2]{BITV}. They then show that the projection map $\mathcal{Z}\mathcal{H}^{(1)}_{\mathfrak{o}} \to \mathcal{Z}\mathcal{H}^{(1)}_{\mathcal{L}}$ is a monoidal equivalence. Therefore in Theorem 5.3.2 of loc.\ cit., the authors identify $D^b_{\mathfrak{C}, \mathcal{L}}(G)$ with the center of $\mathcal{H}^{(1)}_{\mathcal{L}}$. The following proposition is a consequence of this theorem, and it will be an important tool in the present work.

\begin{prop}\label{prop:characterconv}
    There is a well-defined convolution functor
    \begin{align}
         - * - : D_{\mathcal{L}}^b(G/U)\times D^b_{\mathfrak{C}, \mathcal{L}}(G) \to D_{\mathcal{L}}^b(G/U),
    \end{align}
    and for any $Z \in D^b_{\mathfrak{C}, \mathcal{L}}(G)$, the functor $- * Z$ is central: it commutes with convolution by elements of ${}_{\mathcal{L}}\hat{\mathcal{P}}_{\mathcal{L}}$ or ${}_{\mathcal{L}}{\mathcal{P}}_{\mathcal{L}}$ and the isomorphism realizing this commutativity satisfies the corresponding associativity constraints.
\end{prop}

In \Cref{def:charsheaves}, clearly $D^b_{\mathfrak{C},\mathcal{L}}(G) = D^b_{\mathfrak{C},\mathcal{L}'}(G)$ as categories whenever $\mathcal{L}$ and $\mathcal{L}'$ are in the same $W$-orbit $\mathfrak{o}$. However, the convolution action described in \Cref{prop:characterconv} passes through the identification of $D^b_{\mathfrak{C}, \mathcal{L}}(G)$ with the center of $\mathcal{H}^{(1)}_{\mathcal{L}}$, and so by this convention, this action which we will use throughout the present paper depends on $\mathcal{L}$ itself and not just its orbit. To make this identification explicit, we use the following definition.

\begin{defn}
    Let $\mathfrak{hc}_{\mathcal{L}} : D^b_{\mathfrak{C}, \mathcal{L}}(G) \to \mathcal{H}_{\mathcal{L}}^{(1)}$ be the composition of $\mathfrak{hc}$ with the projection $\mathcal{H}_{\mathfrak{o}}^{(1)} \to \mathcal{H}_{\mathcal{L}}^{(1)}$.
\end{defn}

\subsubsection{Two-sided cells and character sheaves}

In this and subsequent sections, we use the notion of two-sided Kazhdan-Lusztig cells in the Weyl group; see e.g. \cite{W} for a clear exposition.

\begin{defn}
    For any $\mathcal{L}$, let $\underline{C}_{\mathcal{L}}$ denote the set of two-sided cells in $W_{\mathcal{L}}^\circ$. For any $w \in W_{\mathcal{L}}^\circ$, we let $\underline{c}_w$ denote the corresponding cell. We let $\underline{c}_{e}$ be the top cell which corresponds to the identity element. There is a well-defined partial order $\leq$ on $\underline{C}_{\mathcal{L}}$ for which $\underline{c}_e$ is maximal.
\end{defn}

Two-sided Kazhdan-Lusztig cells give a filtration on the category $D^b_{\mathfrak{C}, \mathcal{L}}(G)$ (c.f. \cite{CSIV}), and so for each $\underline{c} \in \underline{C}_{\mathcal{L}}$, there are triangulated subcategories $D^b_{\mathfrak{C}, \mathcal{L}}(G)_{\leq \underline{c}}$ and $D^b_{\mathfrak{C}, \mathcal{L}}(G)_{< \underline{c}}$ of $D^b_{\mathfrak{C}, \mathcal{L}}(G)$. We then define $D^b_{\mathfrak{C}, \mathcal{L}}(G)_{\underline{c}}$ as the quotient category $D^b_{\mathfrak{C}, \mathcal{L}}(G)_{\leq \underline{c}}/D^b_{\mathfrak{C}, \mathcal{L}}(G)_{< \underline{c}}$, referring to the unipotent case treated in Section 5 of \cite{BFO} for details. Let $G_{\mathrm{ad}}$ be the adjoint quotient of $G$; the following is a consequence of the classification of irreducible character sheaves in terms of cells given in \cite{CSIV}, c.f. Corollary 5.4 of \cite{BFO}. 
\begin{prop}\label{prop:cellsum} For any $\mathcal{L} \in \mathrm{Ch}(T)$,
    \begin{align}
        K_0(D^b_{\mathfrak{C}, \mathcal{L}}(G_{\mathrm{ad}})) \cong \bigoplus_{\underline{c} \in \underline{C}_{\mathcal{L}}} K_0(D^b_{\mathfrak{C}, \mathcal{L}}(G_{\mathrm{ad}})_{\underline{c}})
    \end{align}
    as vector spaces. Further, for any $\underline{c} \in \underline{C}_{\mathcal{L}}$, the preimage of the subspace
    \begin{align}
        \bigoplus_{\underline{c}' \leq \underline{c}} K_0(D^b_{\mathfrak{C}, \mathcal{L}}{(G_{\mathrm{ad}})}_{\underline{c}'})
    \end{align}
    in $K_0(D^b_{\mathfrak{C}, \mathcal{L}}(G_{\mathrm{ad}}))$ under this isomorphism is a monoidal ideal.
\end{prop}

\subsubsection{The big free-monodromic tilting object and $\mathbb{K}_{\mathcal{L}}$}

In \cite[Section 9.4]{Gouttard}, the author defines free-monodromic tilting sheaves with general monodromy, analogous to the ones appearing in \cite{BY} for the case of unipotent monodromy.
\begin{defn}
    Given $\mathcal{L} \in \mathrm{Ch}(T)$, let $\hat{\mathcal{T}}(w_{0,\mathcal{L}})_{\mathcal{L}}$ be the free-monodromic tilting object in ${}_{\mathcal{L}}\hat{\mathcal{P}}_{\mathcal{L}}$ corresponding to the longest element $w_{0,\mathcal{L}}$ of $W_{\mathcal{L}}^\circ$. In this paper, we will denote it simply by $\hat{\mathcal{T}}_{\mathcal{L}}$ for convenience.
\end{defn}

In \cite{BT}, working in the case of unipotent monodromy, the authors construct from $\hat{\mathcal{T}}(w_0)$ an object they call $\mathbb{K}$, defined as $\mathbb{K} = p^*p_!\hat{\mathcal{T}}(w_0)$ where $p : U\backslash G/U \to (U\backslash G/U)/T$ where $T$ acts on $U\backslash G / U$ by conjugation.

They also define a character sheaf $\Xi$ whose details are explained in 1.4 of \cite{BT} obtained by averaging the derived pushforward of the constant sheaf on the regular locus of the unipotent variety of $G$ to obtain an element of $D^b(G/_{\mathrm{Ad}} G)$. They then show how to define the projection of $\Xi$ to the subcategory of unipotent character sheaves, and that $\mathbb{K}$ is obtained by applying the functor $\mathfrak{hc}$ to the resulting object.
\begin{prop}[(Theorem 1.4.1, \cite{BT})]
    The object $\mathbb{K}$ lies in the image of the functor $\mathfrak{hc}$.
\end{prop}

We now define an anologue for arbitrary monodromy generalizing the unipotent case.
\begin{defn}
    Let $\mathbb{K}_{\mathcal{L}} = p^*p_!\hat{\mathcal{T}}_{\mathcal{L}}$, and let $\mathfrak{o}$ be the $W$-orbit of $\mathcal{L}$. Then by the same argument as in the proof of Theorem 1.4.1 of \cite{BT}, $\mathbb{K}_{\mathcal{L}}$ arises as the image of an object in $D^b_{\mathfrak{C}, \mathcal{L}}(G)$ under the functor $\mathfrak{hc}_{\mathcal{L}}$. (It is the image of the projection of $\Xi$ to the derived category of character sheaves with monodromy in $\mathfrak{o}$.) As a result, we identify $\mathbb{K}_{\mathcal{L}}$ with its preimage under $\mathfrak{hc}_{\mathcal{L}}$ and consider it as an element of $D^b_{\mathfrak{C}, \mathcal{L}}(G)$.
\end{defn}

By \Cref{prop:characterconv}, this means that the functor $- * \mathbb{K}_{\mathcal{L}} : D_{\mathcal{L}}^b(G/U) \to D_{\mathcal{L}}^b(G/U)$ commutes with convolution by elements of ${}_{\mathcal{L}}\hat{\mathcal{P}}_{\mathcal{L}}$.

\subsubsection{Fourier transform and convolution with costandard sheaves}

\begin{prop}\label{prop:ftconv}
    Let $\mathcal{F} \in D^b_{\mathcal{L}}(G/U)$ where $\mathcal{L} \in \mathrm{Ch}(T)$. Then for any $s \in S$,
    \begin{align}
        \Phi_s(\mathcal{F}) & = \begin{cases}\label{eqn:casesofnabla}
            \mathcal{F} * {}_\mathcal{L}\hat{\nabla}(s)_{\mathcal{L}}(\tfrac{1}{2}) & s \in W_{\mathcal{L}}^\circ\\
            \mathcal{F} * {}_{\mathcal{L}}\hat{\nabla}(s)_{s\mathcal{L}} & s \not\in W_{\mathcal{L}}^\circ.
        \end{cases}
    \end{align}
\end{prop}
\begin{proof}
    In the case where $s \in W_{\mathcal{L}}^\circ$, this follows for equivariant sheaves by Proposition 4.3 of \cite{CMFKLCatO}, c.f. 6.3 of \cite{P}, while the second case follows from a similar computation, which is done in the proof of \cite[Proposition 4.1]{CMFSymplectic}. The proof then generalizes to arbitrary extensions of $\mathcal{L}$-equivariant sheaves, and therefore to all monodromic sheaves as in the claim.
\end{proof}

We note that the difference between the two cases in (\ref{eqn:casesofnabla}) boils down to the calculation of $R\Gamma_{c}(T, \mathcal{L}\otimes\mathcal{L}_{\psi})$ in rank $2$, in the case where $\mathcal{L}$ is trivial versus the case where it is nontrivial. As explained in \cite[Applications de la formule des traces aux sommes trigonom{\'e}triques]{DCoh}, this is always concentrated in cohomological degree $1$, but has weight $0$ or $1$ depending on whether $\mathcal{L}$ is trivial, hence the presence of the Tate twist in (\ref{eqn:casesofnabla}). The full computation is explained in more detail in \cite[Proposition 4.1]{CMFSymplectic}.

\section{Coalgebras and dg enhancements}\label{sec:dg}

\subsection{Coalgebras over comonads and Barr-Beck for Kazhdan-Laumon categories}

\subsubsection{The Kazhdan-Laumon category as (co)algebras over a (co)monad}

We recall a result from \cite{BBP} exhibiting the Kazhdan-Laumon category $\mathcal{A}$ as the category of (co)algebras over a certain (co)monad on the underlying category $\mathcal{B} = \mathrm{Perv}(G/U)^{\oplus W}$.

\begin{defn}
    Let $\Psi^\circ : \mathcal{B} \to \mathcal{B}$ be the endofunctor defined for any $A = (A_w)_{w \in W} \in \mathcal{B}$ by
    \begin{align}\label{eqn:comonad}
        (\Psi^\circ A)_w & = \oplus_{y \in W} \Psi_{wy^{-1}}^\circ A_y.
    \end{align}
    It has right-derived functor $\Psi = R\Psi^\circ$ given by
    \begin{align}
        (\Psi A)_w & = \oplus_{y \in W} \Psi_{wy^{-1}} A_y.
    \end{align}
    Similarly, we define $\Phi^\circ : \mathcal{B} \to \mathcal{B}$ by 
    \begin{align}
        (\Phi^\circ A)_w & = \oplus_{y \in W} \Phi_{wy^{-1}}^\circ A_y,
    \end{align}
    with left-derived functor $\Phi = L\Phi^\circ$ given by
    \begin{align}
        (\Phi A)_w & = \oplus_{y \in W} \Phi_{wy^{-1}} A_y.
    \end{align}
\end{defn}

\begin{theorem}[\cite{BBP}]
    The Kazhdan-Laumon category $\mathcal{A}$ is equivalent to the category of coalgebras over the left-exact comonad $\Psi^\circ : \mathcal{B} \to \mathcal{B}$. Dually, $\mathcal{A}$ is equivalent to the category of algebras over the right-exact monad $\Phi^\circ : \mathcal{B} \to \mathcal{B}$. 
\end{theorem}

\subsubsection{dg enhancements of derived categories}

Given an abelian category $\mathcal{C}$, let $C_{\mathrm{dg}}(\mathcal{C})$ be the dg category with objects complexes of sheaves and morphisms the usual complexes of maps between complexes. We define the dg derived category $D_{\mathrm{dg}}(\mathcal{C})$ to be the dg quotient of $C_{\mathrm{dg}}(\mathcal{C})$ by the full subcategory of acyclic objects; its homotopy category is the usual derived category $D(\mathcal{C})$.

The bounded dg derived category $D^b_{\mathrm{dg}}(\mathcal{C})$ is defined to be the full dg subcategory of $D_{\mathrm{dg}}(\mathcal{C})$ consisting of objects which project to $D^b(\mathcal{C})$ when passing to the homotopy category.

\subsubsection{Barr-Beck-Lurie for dg categories applied to the Kazhdan-Laumon category}

We first state the following general result, which follows from the Barr-Beck-Lurie monadicity theorem.

\begin{prop}\label{prop:barrbeck}
    Suppose $\mathcal{C}$ is a Grothendieck abelian category, and let $T = F\circ F^L$ be the monad on $\mathcal{C}$ arising from a monadic adjunction $(F^L, F)$. Let $\mathcal{C}^T$ be the category of algebras over the monad $T$.

    Now suppose also that $F : \mathcal{C}^T \to \mathcal{C}$ is exact and also admits a right adjoint $F^R$. Suppose also that the left- and right-derived functors $LF^L$ and $RF^R$ give adjunctions $(LF^L, F)$ and $(F, RF^R)$ of functors between the dg derived categories $D_{\mathrm{dg}}(\mathcal{C}^T)$ and $D_{\mathrm{dg}}(\mathcal{C})$. Let $\tilde{T} = F \circ LF^L$ be the resulting dg monad, and $D_{\mathrm{dg}}(\mathcal{C})^{\tilde{T}}$ the dg category of algebras over the monad $\tilde{T}$.
    
    If $F : D_{\mathrm{dg}}(\mathcal{C}^T) \to D_{\mathrm{dg}}(\mathcal{C})$ is conservative, then there is a canonical equivalence of dg categories
    \begin{align}
        \tilde{F} : D_{\mathrm{dg}}(\mathcal{C}^T) \to D_{\mathrm{dg}}(\mathcal{C})^{\tilde{T}}
    \end{align}
\end{prop}
\begin{proof}
    Since $F$ and $LF^L$ are both left adjoints, each preserves small colimits, as does $\tilde{T}$. Because it is also conservative, the functor $F : D_{\mathrm{dg}}(\mathcal{C}^T) \to D_{\mathrm{dg}}(\mathcal{C})$ then satisfies the conditions of the Barr-Beck-Lurie theorem \cite[Theorem 4.7.3.5]{ha} for dg categories. A simplification of this theorem in the case where $F$ and $F^L$ both preserve small colimits is explained in \cite[Section 2.2]{gunningham}, and this assumption holds in the present case.
\end{proof}

Note that the monad $\Phi^\circ : \mathcal{B} \to \mathcal{B}$ can be enhanced to a monad $\Phi = L\Phi^\circ : D_{\mathrm{dg}}^b(\mathcal{B}) \to D_{\mathrm{dg}}^b(\mathcal{B})$; i.e. the functor $\Phi$ has a dg enhancement (since it arises from the functors $\Phi_w$ which are themselves defined using the six-functor formalism). We can then consider the dg category $D_{\mathrm{dg}}^b(\mathcal{B})^{\Phi}$ of algebras over the monad $F \circ \Phi$ or the dg category $D_{\mathrm{dg}}^b(\mathcal{B})_{\Psi}$ of coalgebras over the comonad $F \circ \Psi$. The following is an application of \Cref{prop:barrbeck}.

\begin{prop}\label{prop:dgequiv}
    The dg categories $D_{\mathrm{dg}}^b(\mathcal{B}^{\Phi^\circ})$ and $D_{\mathrm{dg}}^b(\mathcal{B})^{\Phi}$ are equivalent.
\end{prop}
\begin{proof}
    We check the conditions of \Cref{prop:barrbeck} for $\mathcal{C} = \mathcal{B}$ and $F$ being the forgetful functor $\mathcal{B}^{\Phi^\circ} \to \mathcal{B}$ or its dg version $D_{\mathrm{dg}}(\mathcal{B}^{\Phi^\circ}) \to D_{\mathrm{dg}}(\mathcal{B})$, and $T = F\circ \Phi^\circ$. It is clear that $F$ is exact; we claim that $F : D_{\mathrm{dg}}(\mathcal{B}^{\Phi^\circ}) \to D_{\mathrm{dg}}(\mathcal{B})$ has a left adjoint given by the free algebra functor $\mathcal{B} \to \mathcal{B}^{\Phi^\circ}$. 

    Indeed, we check this adjunction explicitly: this free algebra functor sends $(B_w)_{w \in W} \in D^b_{\mathrm{dg}}(\mathcal{B})$ to the direct sum $\oplus_{w \in W} Lj_{w!}^\circ(A_w)$ where $Lj_{w!}^\circ = j_{w!}$ is the left adjoint to $j_w^*$ by \Cref{prop:adjunctionshriek}. It is then clear from the adjunction $(j_{w!}, j_{w}^*)$ that the forgetful functor is its left adjoint. Now notice that in a similar way, we observe that the forgetful functor $F$ has a right adjoint: it follows from the opposite adjunction $(j_w^*, j_{w*})$ that the functor sending any $(B_w)_{w \in W} \in D^b_{\mathrm{dg}}(\mathcal{B})$ to $\oplus_{w \in W} j_{w*}(B_w)$ (where again $j_{w*} = Rj_{w*}^\circ$) is right-adjoint to $F$.

    We are now in the setup of \Cref{prop:barrbeck}, and both the abelian and dg versions of the functor $F$ are clearly conservative. So applying this proposition, we get a natural equivalence $D_{\mathrm{dg}}(\mathcal{B}^{\Phi^\circ}) \to D_{\mathrm{dg}}(\mathcal{B})^{\Phi}$. To conclude, we then note that this then restricts to an equivalence of bounded derived categories $D_{\mathrm{dg}}^b(\mathcal{B}^{\Phi^\circ}) \to D_{\mathrm{dg}}^b(\mathcal{B})^{\Phi}$ since cohomology in both categories is computed in the underlying category $D(\mathcal{B})$ after forgetting the algebra structure.
\end{proof}

Dually, in terms of the comonad $\Psi$, we have that the dg categories $D_{\mathrm{dg}}^b(\mathcal{B}_{\Psi^\circ})$ and $D_{\mathrm{dg}}^b(\mathcal{B})_{\Psi}$ are equivalent.

\subsection{The center of the monodromic Hecke category}

\subsubsection{Action of the center on the Kazhdan-Laumon dg-category}

For any $\mathcal{L} \in \mathrm{Ch}(T)$, let $\mathcal{B}_{\mathcal{L}} = \oplus_{w \in W}\mathrm{Perv}_{w\mathcal{L}}(G/U)$.

\begin{defn}
    Given $Z \in D^b_{\mathfrak{C}, \mathcal{L}}(G)$, we consider $Z$ as a functor on the direct sum of categories $\oplus_{w\in W} \mathrm{Perv}_{w\mathcal{L}}(G/U)$ as follows. Given any object $A = (A_w)$ in $\oplus_{w \in W} \mathrm{Perv}_{w\mathcal{L}}(G/U)$, let $Z(A) = A * Z \in \oplus_{w \in W} \mathrm{Perv}_{w\mathcal{L}}(G/U)$ be defined such that $(A * Z)_w = A_w * {}_{w\mathcal{L}}\hat{\Delta}(w)_{\mathcal{L}} * Z * {}_{\mathcal{L}}\hat{\nabla}(w^{-1})_{w\mathcal{L}}$.
\end{defn}

By \Cref{prop:characterconv}, for any $B \in {}_{\mathcal{L}}\hat{\mathcal{P}}_{\mathcal{L}}$, the functors $- * Z * B$ and $- * B * Z$ on $\mathrm{Perv}_{\mathcal{L}}(G/U)$ are naturally isomorphic. The following lemma is a consequence of this fact.
\begin{lemma}\label{lem:naturaliso}
    For any $Z \in D^b_{\mathfrak{C}, \mathcal{L}}(G)$, there is a natural isomorphism $Z \circ \Psi \cong \Psi \circ Z$ of endofunctors of $\mathcal{A}^{\mathcal{L}}$.
\end{lemma}
\begin{proof} By \Cref{prop:ftconv}, for any $w \in W$, there exists some Tate twist $(d)$ such that $((Z \circ \Psi)(A))_w(d)$ is isomorphic to
{\allowdisplaybreaks
    \begin{align}
        & \quad \oplus_{y \in W}  (\Psi_{wy^{-1}} A_y) * {}_{w\mathcal{L}}\hat{\Delta} (w)_{\mathcal{L}} * Z * {}_{\mathcal{L}}\hat{\nabla}(w^{-1})_{w\mathcal{L}}\\
        & \cong  A_y * {}_{y\mathcal{L}}\hat{\nabla}(yw^{-1})_{w\mathcal{L}} * {}_{w\mathcal{L}}\hat{\Delta}(w)_{\mathcal{L}} * Z * {}_{\mathcal{L}}\hat{\nabla}(w^{-1})_{w\mathcal{L}}\\
        & \cong A_y * {}_{y\mathcal{L}}\hat{\Delta}(y)_{\mathcal{L}} * {}_{\mathcal{L}}\hat{\nabla}(y^{-1}){}_{y\mathcal{L}} * {}_{y\mathcal{L}}\hat{\nabla}(yw^{-1})_{w\mathcal{L}} * {}_{w\mathcal{L}}\hat{\Delta}(w)_{\mathcal{L}} * Z * {}_{\mathcal{L}}\hat{\nabla}(w^{-1})_{w\mathcal{L}} \\
        & \cong A_y * {}_{y\mathcal{L}}\hat{\Delta}(y){}_{\mathcal{L}} * Z * {}_{\mathcal{L}}\hat{\nabla}(y^{-1}){}_{y\mathcal{L}} * {}_{y\mathcal{L}}\hat{\nabla}(yw^{-1}){}_{w\mathcal{L}} * {}_{w\mathcal{L}}\hat{\Delta}(w){}_{\mathcal{L}} * {}_{\mathcal{L}}\hat{\nabla}(w^{-1}){}_{w\mathcal{L}}\\
        & \cong A_y * {}_{y\mathcal{L}}\hat{\Delta}(y){}_{\mathcal{L}} * Z * {}_{\mathcal{L}}\hat{\nabla}(y^{-1}){}_{y\mathcal{L}} * {}_{y\mathcal{L}}\hat{\nabla}(yw^{-1}){}_{w\mathcal{L}} \\
        & \cong \oplus_{y \in W} \Psi_{wy^{-1}}(Z(A))(d)\\
        & \cong ((\Psi \circ Z)(A))_w(d),
    \end{align}}
    and these isomorphisms across all $w \in W$ assemble together to give a natural isomorphism $Z \circ \Psi \cong \Psi \circ Z$. A conceptual explanation for the existence of this isomorphism appears in Lemma 11.12 of \cite{lusztig_yun_2020}.
\end{proof}

In the following, let $\pi_{\mathrm{ad}} : G \to G_{\mathrm{ad}}$ be the natural map from $G$ to its adjoint quotient. 
\begin{prop}\label{prop:dgacts}
    There is an action of the monoidal category $D^b_{\mathfrak{C},\mathcal{L}}(G_{\mathrm{ad}})$ on the category $D^b(\mathcal{A}^{\mathcal{L}})$, which we denote by convolution on the right.
\end{prop}
\begin{proof}
    Given an element $Z \in D^b_{\mathfrak{C},\mathcal{L}}(G_{\mathrm{ad}})$, we first pull back along $\pi_{\mathrm{ad}}$ to obtain an element $\tilde{Z} \in D^b_{\mathfrak{C},\mathcal{L}}(G)$. We can then act by the usual convolution $- * \tilde{Z}$ described in \Cref{prop:characterconv}. We note that $D^b_{\mathfrak{C},\mathcal{L}}(G_{\mathrm{ad}})$ has a dg enhancement discussed for the unipotent case in \cite{BZN} but analogous in general. Further, the functors we use are built from compositions of those occuring in the six-functor formalism and therefore each has a dg enhancement as well.

    For any $B \in D^b_{\mathrm{dg}}(\mathcal{B}_\mathcal{L})_{\Psi}$ and $Z \in  D^b_{\mathfrak{C},\mathcal{L}}(G_{\mathrm{ad}})$, we define the new object $B * Z$ by $F(B) * Z$ as a tuple of elements of $\mathrm{Perv}_{\mathcal{L}}(G/U)$ where $F$ is the forgetful functor. We then define the coalgebra structure map $B * Z \to \Psi (B * Z)$ by the composition
    \begin{align}
        B * Z & \to \Psi(B) * Z \cong \Psi(B * Z)
    \end{align}
    where the first map comes from the coalgebra structure on $B$ and the subsequent isomorphism is the natural isomorphism described in Lemma \ref{lem:naturaliso}. One can then check that this defines a coalgebra structure on $B * Z$.
\end{proof}

We note that we get a similar action of $D^b_{\mathfrak{C},\mathcal{L}}(G_{\mathrm{ad}})$ on $D^b(\AwL)$ since any object $Z$ here is a sheaf on the variety $G_{\mathrm{ad}}$ and therefore has a natural isomorphism $\mathrm{Fr}^* Z \to Z$; this means $- * Z$ commutes with $\mathrm{Fr}^*$, and therefore the action in \ref{prop:dgacts} also gives an action in the setting of $w$-twisted Weil sheaves.

\section{Polishchuk's canonical complex}\label{sec:canonical}

\subsection{Parabolic Kazhdan-Laumon categories}
\subsubsection{Definition of the categories $\mathcal{A}_{w, \mathbb{F}_q}^I$}
We will now define certain parabolic analogues of the Kazhdan-Laumon category which correspond to subsets $I \subset S$. We first recall the definition of the categories $\mathcal{A}_{W_I}$ from \cite{P}, build off of this definition to define the categories $\mathcal{A}^I$ which we will need in the present work, and finally define $\mathcal{A}^I_{w, \mathbb{F}_q}$ for any $w \in W$.

\begin{defn}[(\cite{P}, Section 7)]
    For any $I \subset S$, the category $\mathcal{A}_{W_I}$ is the category of tuples of elements of $\mathrm{Perv}(G/U)$ indexed by $W_I$ with morphisms and compatibilities as in \Cref{def:kldef} but only for $w \in W_I$, $s \in I$.

    For any $J \subset I \subset S$, and any right coset $W_Jx$ of $W_J$ in $W_I$, there is a restriction functor $j_{W_Jx}^{W_I*} : \mathcal{A}_{W_I} \to \mathcal{A}_{W_J}$ remembering only the tuple elements and morphisms for $w \in W_Jx$, $s \in J$. When $I = S$, we write only $j_{W_Jx}^*$, and in this case we omit superscripts similarly for all functors introduced in this section.
\end{defn}

\begin{prop}[(\cite{P}, Proposition 7.1.1)]\label{prop:polparind}
    For any $J \subset I \subset S$, and any $W_Jx$, the functor $j_{W_Jx}^{W_I*}$ admits a left adjoint $j_{W_Jx!}^{W_I}$.
\end{prop} 

\begin{defn}
    For $I \subset S$, let $\mathcal{A}^I$ be the category whose objects are tuples $(A_w)_{w \in W}$ equipped with morphisms $\Phi_y^\circ A_w \to A_{yw}$ for any $w \in W$, $y \in W_I$ satisfying the same conditions as in \Cref{def:kldef} but only for $y \in W_I$. Morphisms in $\mathcal{A}^I$ are morphisms in $\mathrm{Perv}(G/U)^{\oplus W}$ satisfying the compatibilities in \Cref{def:kldef} but only for $y \in W_I$.

    Alternatively, $\mathcal{A}^I = \oplus_{W_I \backslash W} \mathcal{A}_{W_I}$ with a reindexing by $W$ on the tuples in this category.
\end{defn}

Just like $\mathcal{A}$, the category $\mathcal{A}^I$ admits an action of $W$ by functors $\{\mathcal{F}_w\}_{w \in W}$ which are defined by $\mathcal{F}_w((A_y)_{y \in W}) = (A_{yw})_{w \in W}$. We define the category $\mathcal{A}_{w, \mathbb{F}_q}^I$ as in \Cref{def:aw} but with $\mathcal{A}$ replaced by $\mathcal{A}^I$.
\subsubsection{Adjunctions between these categories}
For any $J \subset I$, there is an obvious restriction functor $j_{I,J}^* : \mathcal{A}^I \to \mathcal{A}^J$ which is the identity on objects but which remembers only the morphisms $\Phi_y^\circ A_w \to A_{yw}$ for $y \in W_J$. As in \Cref{prop:adjunctionshriek}, there is an analogous derived version of this morphism and the following adjunction.
\begin{prop}
    For any $J \subset I$, the functor $j_{I,J}^*$ admits a left adjoint $j_{I,J!}$. For any $A \in \mathcal{A}^J$, $j_{I,J!}^\circ(A)$ is the direct sum
\begin{align}\label{eqn:dirsumadj}
    \bigoplus_{x \in W_J\backslash W_I} j^{W_I\circ}_{W_Jx!}((A_w)_{w \in W_Jx}),
\end{align}
where $(A_w)_{w \in W_Jx}$ is considered as an element of $\mathcal{A}_{W_J}$.
Further, the adjoint pair $(j_{I,J!}^\circ, j_{I,J}^*)$ gives also an adjunction between $\mathcal{A}^I_{w, \mathbb{F}_q}$ and $\mathcal{A}^J_{w, \mathbb{F}_q}$. 
\end{prop} 

\begin{proof}
The direct sum in (\ref{eqn:dirsumadj}) has the natural structure of an object of $\mathcal{A}^I$, as each summand object $j_{W_Jx!}^{W_I \circ}(A_w)_{w \in W_J x}$ has such a structure by Proposition \ref{prop:polparind}. We then note that for any such $A \in \mathcal{A}^J$ and any $B \in \mathcal{A}^I$,
{\allowdisplaybreaks
\begin{align*}
& \quad \mathrm{Hom}_{\mathcal{A}^I}\left(\bigoplus_{x \in W_J \backslash W_I} j_{W_Jx!}^{I\circ}((A_w)_{w \in W_Jx}), B\right)\\
& \cong \bigoplus_{x \in W_J \backslash W_I} \mathrm{Hom}_{\mathcal{A}^I}(j_{W_J x!}^\circ((A_w)_{w \in W_J x}), B)\\
& \cong \bigoplus_{x \in W_J \backslash W_I} \mathrm{Hom}_{\mathcal{A}_{W_J x}}((A_w)_{w \in W_J x}, j_{W_J x}^* B)\\
& \cong \mathrm{Hom}_{\mathcal{A}^J}((A_w)_{w \in W}, \oplus_{x \in W_J \backslash W_I}j_{W_J x}^*B)\\
& = \mathrm{Hom}_{\mathcal{A}^J}((A_w)_{w \in W}, j_{I,J}^*B).
\end{align*}}
To see that the adjoint pair $(j_{I,J!}^\circ, j_{I,J}^*)$ gives also an adjunction between $\mathcal{A}^I_{w, \mathbb{F}_q}$ and $\mathcal{A}^J_{w, \mathbb{F}_q}$, we note that the morphisms on both sides of the equation above which are compatible with the morphism $\psi_A : \mathcal{F}_w\mathrm{Fr}^* A \to A$ are preserved by these isomorphisms.
\end{proof}

\subsection{Polishchuk's complex for $\mathcal{A}_{w, \mathbb{F}_q}$}
\subsubsection{Polishchuk's canonical complex\label{subsubsec:pcc}}

For a fixed $A \in \mathcal{A}$ and a choice of $J \subset S$, Polishchuk writes $A(J) = j_{S - J!}j_{S - J}^*A$. In 7.1 of \cite{P}, Polishchuk explains that for any $A \in \mathcal{A}$, adjunction of parabolic pushforward and pullback functors $(j_{W_Ix!}, j_{W_Ix}^*)$ gives a canonical morphism $A(J) \to A(J')$ whenever $J \subset J'$. He then defines the complex $C_\bullet(A)$ as
\[\begin{tikzcd}
    C_{n-1} = A(S) \arrow[r] & \dots \arrow[r] & C_1 = \bigoplus_{|J| = 2} A(J) \arrow[r] & C_0 = \bigoplus_{|J| = 1} A(J),
\end{tikzcd}\]
where $n = |S|$.

Recalling the morphism $\iota : \mathcal{A} \to \mathcal{A}$ sending $(A_w)_{w \in W}$ to $(\Phi_{w_0}^\circ A_{w_0w})_{w \in W}$, he describes natural morphisms $C_0 \to A$ and $\iota(A) \to C_{n-1}$, and shows that
\begin{align}
    H_i(C_\bullet(A)) & = \begin{cases}
        A & i = 0,\\
        0 & i \neq 0, n - 1,\\
        \iota(A) & i = n-1.
    \end{cases}
\end{align}
He then describes a $2n$-term complex $\tilde{C}_\bullet$ formed from attaching $C_\bullet(\iota A)$ to $C_\bullet(A)$ via the maps $C_0(\iota A) \to \iota A \to C_{n-1}(A)$ with the property that
\begin{align}
    H_i(\tilde{C}_{\bullet}(A)) & = \begin{cases}
        A & i = 0,\\
        0 & i \neq 0, 2n - 1,\\
        \iota^2(A) & i = 2n - 1.
    \end{cases}
\end{align}

\subsubsection{Compatibility of the complex with $w$-twisted structure}

\begin{prop}
    If $A \in \mathcal{A}_{w, \mathbb{F}_q}$, then the complex $C_\bullet(A)$ is compatible with the $w$-twisted Weil structures, i.e. is a complex of objects in $\mathcal{A}_{w, \mathbb{F}_q}$.
\end{prop}

\begin{proof}
    For any $k$, components of the map $C_k(A) \to C_{k-1}(A)$ are each maps of the form
    \begin{align}\label{eqn:polimap}
        j_{S-J!}j_{S-J}^*A \to j_{S-J'!}j_{S-J'}^*A
    \end{align}
    where $J' \subset J \subset S$ are such that $|J| = k$, $|J'| = k - 1$, which we now describe. By adjunction the data of such a map is equivalent to a map
    \begin{align}\label{eqn:poliadj}
        j_{S-J}^*A \to j_{S-J}^*j_{S-J'!}j_{S-J'}^*A,
    \end{align}
    and compatibility with the $w$-twisted structure is preserved under this adjunction. By the definition of $j_{S-J'!}$, we have that
    \begin{align}
        j_{S-J}^*j_{S-J'!}j_{S-J'}^*A & \cong \bigoplus_{x \in W_{S-J}\backslash W} j_{S-J}^*j_{W_{S-J}x!}j_{S-J'}^*A,
    \end{align}
    and one can check that the map in (\ref{eqn:polimap}) appearing in the definition of Polishchuk's complex in \cite{P} comes in (\ref{eqn:poliadj}) from a natural injection in $\mathcal{A}^{S-J}$ from $j_{S-J}^*A$ into this direct sum defined by sending the $y$th tuple entry to the $y$th tuple entry in the direct summand corresponding to the unique $x$ for which $y \in W_{S-J}x$. It is straightforward to check that this injection preserves $w$-twisted Weil structures on both sides coming from the $w$-twisted Weil structure on $A$, and therefore the map in (\ref{eqn:polimap}) is a map in $\mathcal{A}_{w, \mathbb{F}_q}$. 
\end{proof}

\subsubsection{Parabolic canonical complexes} We remark how that the content appearing in \ref{subsubsec:pcc} can be generalized to provide complexes in $\mathcal{A}^I_{w, \mathbb{F}_q}$ for any $I \subset S$ with $|I| = k$ and any $w \in W$. Namely, if we fix $I \subset S$, $w \in W$, and $A \in \mathcal{A}_{w, \mathbb{F}_q}^I$, we let $A^I(J) = j_{S-J!}^Ij_{S-J}^{I*}A$ whenever $J \supset S - I$, and we define the complex $C^I_\bullet(A)$ by
\[\begin{tikzcd}
    A^I(S) \arrow[r] & \dots \arrow[r]  & {\bigoplus_{|J|=n-k + 2}}A^I(J) \arrow[r] & \bigoplus_{|J| = n-k + 1}A^I(J),
\end{tikzcd}\]
indexed such that $C_{k-1}$ is the first term and $C_0$ is the last term in the above. The results in \Cref{subsubsec:pcc} then still hold, giving a version of the canonical complex associated to an object $A \in \mathcal{A}_{w, \mathbb{F}_q}^I$. 

\subsubsection{Equations in the Grothendieck group} The following is a consequence of the fact that the full twist is central in the braid group, along with the fact that the symplectic Fourier transforms $\Phi_w$ form a braid action.
\begin{lemma}\label{lem:iotacommutes}
    For any $J \subset I \subset S$ there is a natural isomorphism
    \begin{align}
        j_{J!}^I\circ \iota^2 \cong \iota^2\circ j_{J!}^I
    \end{align}
    of functors from $D^b(\mathcal{A}^J)$ to $D^b(\mathcal{A}^I)$.
\end{lemma}

\begin{theorem}\label{thm:i2min1}
    For any $w \in W$ and $A \in \mathcal{A}_{w, \mathbb{F}_q}$ the element
    \begin{align}
        (\iota^2 - 1)^n[A] \in K_0(\mathcal{A}_{w, \mathbb{F}_q})
    \end{align}
    lies in $\vfp$, for $n = |S|$.
\end{theorem}

\begin{proof}
    By the derived category version of the canonical complex construction in combination with the observation about its homology provided in \ref{subsubsec:pcc}, for any $A \in \mathcal{A}_{w, \mathbb{F}_q}$ we have the following equation in $K_0(\mathcal{A}_{w, \mathbb{F}_q})$:
    \begin{align}
        [\iota A] + (-1)^{|I| - 1}[A] & = \sum_{\substack{I'\\ I \subset I' \subsetneq S}} (-1)^{|I'|} [j_{J!}^Ij_{J}^{I*} A],
    \end{align}
c.f. the proof of Theorem 11.5.1 in \cite{P} where the analogous equation is used in the case where $w = e$, $I = S$. By the ``doubled" canonical complex $\tilde{C}_{\bullet}(A)$ and the description of its homology in \ref{subsubsec:pcc}, this means that $[\iota^2 A] - [A]$ is a linear combination of elements lying in the image of the functors $j_{J!}^Ij_{J}^{I*}$ for $J \subsetneq I$.

Now by induction on the $|I|$ appearing in the equation above, it follows from Lemma \ref{lem:iotacommutes} that $(\iota^2 - 1)^n[A]$ is a linear combination of elements lying in the image of the functors $j_{\varnothing!}j_{\varnothing}^*$. We have that for any $B \in \mathcal{A}^\varnothing_{w, \mathbb{F}_q}$,
\begin{align}
    j_{\varnothing!}j_{\varnothing}^* = \oplus_{y \in W} j_{y!}B_y,
\end{align}
and each of these direct summands has finite cohomological dimension by \Cref{prop:adjunction}, completing the proof of the theorem.
\end{proof}

\section{Central objects and the full twist}\label{sec:central}

\subsection{Cells, the big tilting object, and the full twist}

\subsubsection{The full twist}

\begin{defn}
    Let $\mathcal{L} \in \mathrm{Ch}(T)$. Then we define the element
    \begin{align}
        \mathrm{FT}_{\mathcal{L}} & = {}_{\mathcal{L}}\hat{\nabla}(w_0)_{w_0\mathcal{L}} * {}_{w_0\mathcal{L}}\hat{\nabla}(w_0)_{\mathcal{L}}
    \end{align}
    of ${}_{\mathcal{L}}\hat{\mathcal{P}}_{\mathcal{L}}$. The object $\mathrm{FT}_{\mathcal{L}}$ admits a central structure and arises comes from a character sheaf in $D^b_{\mathfrak{C}, \mathcal{L}}(G_{\mathrm{ad}})$ via pullback composed with $\mathfrak{hc}_{\mathcal{L}}$, see \cite{BT} for an explicit description of this character sheaf in the unipotent case, whose proof can be adapted similarly for arbitrary monodromy. As with $\mathbb{K}_{\mathcal{L}}$, we will identify $\mathrm{FT}_{\mathcal{L}}$ with its underlying character sheaf in $D^b_{\mathfrak{C}, \mathcal{L}}(G_{\mathrm{ad}})$, and by $- * \mathrm{FT}_{\mathcal{L}}$ we will denote the action as described in \Cref{prop:dgacts}.
\end{defn}

We note that to properly extend the result in \cite{BT} to the case of non-unipotent monodromy, one needs an analogue of \cite[Theorem 5.6.4]{BT} (the statement that the Harish-Chandra transform composed with the long intertwining functor, i.e.\ the $*$-version of the Radon transform, commutes with Verdier duality) adapted to the $\ell$-adic setting. Before loc.\ cit., this result appeared for unipotent monodromy in \cite[Corollary 7.9]{CYD}. As the authors of \cite{BITV} note, the methods of \cite{CYD} can be adapted to treat the case of general monodromy. They also suggest that a new uniform proof will appear in upcoming work. Finally, note also that \cite[Corollary 3.4]{BFO} treats the case of general monodromy in the characteristic $0$ setting. 

By the definition of $\mathrm{FT}_{\mathcal{L}}$ combined with \Cref{prop:ftconv}, we obtain the following.
\begin{lemma}\label{lem:i2ft}
    The functors $- * \mathrm{FT}_{\mathcal{L}}$ and $\iota^2$ on $D^b(\AwL)$ are naturally isomorphic.
\end{lemma}

\subsubsection{The action of the full twist on cells}

\begin{prop}\label{prop:celleigvals}
For any $\underline{c} \in \underline{C}_{\mathcal{L}}$, there exists an integer $d_{\mathcal{L}}(\underline{c})$ between $0$ and $2\ell(w_0)$ such that for any $a \in K_0(D^b_{\mathfrak{C}, \mathcal{L}}(G)_{\underline{c}})$,
\begin{align}
    [\mathrm{FT}_{\mathcal{L}}] * a = v^{d_{\mathcal{L}}(\underline{c})}a
\end{align}
\end{prop}
\begin{proof}
    We follow the same argument as in Proposition 4.1 and Remark 4.2 of \cite{BFO}; in other words, we begin with the observation that $[\mathrm{FT}_{\mathcal{L}}]$ acts trivially on the non-graded version of the Grothendieck group $K_0(D^b_{\mathfrak{C}, \mathcal{L}}(G)_{\underline{c}})$. Continuing to follow the argument of loc. cit., we then know that for any object $A$ in the heart of $D^b_{\mathfrak{C}, \mathcal{L}}(G)_{\underline{c}}$, the object $\mathrm{FT}_{\mathcal{L}} * A \in D^b_{\mathfrak{C}, \mathcal{L}}(G)_{\underline{c}}$ is perverse up to shift, and furthermore has the property that $[\mathrm{FT}_{\mathcal{L}} * A] = v^{d}[A]$ for some $d$.

    To compute the value of $d$, we can pass back along the Harish-Chandra transform and work in the category $\mathcal{H}_{\mathcal{L}}^{(1)}$. The Grothendieck ring $K_0(\mathcal{H}_{\mathcal{L}}^{(1)})$ is the monodromic Hecke algebra $\mathcal{H}_{\mathcal{L}}$. By \cite{lusztig_yun_2020}, this is isomorphic to the usual Hecke algebra associated to the group $W_{\mathcal{L}}^\circ \subset W$, with $\mathrm{FT}_{\mathcal{L}}$ being identified with the usual full-twist $\tilde{T}_{w_{0,\mathcal{L}}}^2$. By \cite[5.12.2]{Lbook}, the full twist in the usual Hecke algebra acts on the cell subquotient module of the Hecke algebra corresponding to a cell $\underline{c}$ by the scalar $v^{d(\underline{c})}$, where $d(\underline{c})$ is described in loc. cit. Passing this fact back along the monodromic-equivariant isomorphism from \cite{lusztig_yun_2020}, the result follows.
\end{proof}

\begin{defn}
    For any two-sided cell $\underline{c}$, let $d_{\mathcal{L}}(\underline{c})$ be the integer between $0$ and $2\ell(w_0)$ for which the equation in \Cref{prop:celleigvals} holds.
\end{defn}

\subsubsection{$\mathbb{K}_{\mathcal{L}}$ in the top cell subquotient}

\begin{defn}
    Let $K_0(D^b(\mathcal{A}_{w, \mathbb{F}_q}^{\mathcal{L}}))_{< \underline{c}_e}$ be the submodule of $K_0(D^b(\AwL))$ spanned by the image of the ideal $K_0(D^b_{\mathfrak{C},\mathcal{L}}(G_{\mathrm{ad}})_{< \underline{c}_{e}})$ under the action described in \Cref{prop:dgacts}.
\end{defn}

\begin{defn}\label{def:qv}
    For any $\mathcal{L} \in \mathrm{Ch}(T)$, let $q_{\mathcal{L}}(v)$ be the Poincar{\'e} polynomial
    \[q_{\mathcal{L}}(v) = \sum_{w \in W} (-v^2)^{\ell(w)}\]
    of the group $W_{\mathcal{L}}^\circ$. Note by the Chevalley-Solomon formula that $q_{\mathcal{L}}(v)$ can be expressed as the product of some linear factors each of which is a factor of $(v^{2i} - 1)$ for some $1 \leq i \leq \ell(w_0)$.\footnote{The variable $u$ used throughout \cite{P} is replaced in the present work by $v^2$.}
\end{defn}

\begin{lemma}\label{lem:yuntilting}
    The multiplicity with grading of the irreducible object $\mathrm{IC}(e)_{\mathcal{L}}$ in the Jordan-Holder decomposition of $\hat{\mathcal{T}}_{\mathcal{L}}$ is $q_{\mathcal{L}}(v)$.
\end{lemma}
\begin{proof}
    In \cite{Z}, Yun computes the $\mathbb{Z}[v,v^{-1}]$-graded multiplicity of any standard object $\Delta(y)$ in a filtration of $T(w)$, where $w, y \in W$ and $\Delta(y)$ and $T(w)$ are standard and tilting objects respectively in the usual Hecke category. The main equivalence of categories in \cite{lusztig_yun_2020} allows us to extend these results to the monodromic Hecke category by replacing $W$ with $W_{\mathcal{L}}^\circ$, whose combinatorics in terms of tilting, standard, and irreducible objects matches exactly the combinatorics of the completed category ${}_{\mathcal{L}}\hat{\mathcal{P}}_\mathcal{L}$ (c.f. 9.3.3 of \cite{Gouttard} for an explicit description of the standard filtration on a tilting object in the monodromic setting).

    Combining Theorem 5.3.1 of \cite{Z} with the expression of standard objects in terms of irreducible objects via inverse Kazhdan-Lusztig polynomials, we compute that the multiplicity of $\hat{\mathrm{IC}(e)}_{\mathcal{L}}$ is exactly
    \begin{align}
        \sum_{w \in W_{\mathcal{L}}^\circ} (-v^2)^{\ell(w_0) - \ell(w)} = \sum_{w \in W_{\mathcal{L}}^\circ} (-v^2)^{\ell(w)} = q_{\mathcal{L}}(v).
    \end{align}
\end{proof}

For the next proposition, we recall the definition of the character sheaves $\varepsilon_{n,\mathcal{L}}$ from \Cref{subsubsec:free-monodromic}.
\begin{prop}\label{prop:tiltingid}
For any $n$, the element
\begin{align}\label{eqn:tiltid}
    [\varepsilon_{n,\mathcal{L}} * \mathbb{K}_{\mathcal{L}}] - (v^2 - 1)^{\mathrm{rank}(T)}q_{\mathcal{L}}(v)[\varepsilon_{n,\mathcal{L}}]
\end{align}
    of $K_0(D^b_{\mathfrak{C}}(G_{\mathrm{ad}}))$ lies in the subspace $K_0(D^b_{\mathfrak{C},< \underline{c}_{e}}(G_{\mathrm{ad}}))$.
\end{prop}

\begin{proof}
    First note that $K_0(D^b_{\mathfrak{C},\mathcal{L}}(G_{\mathrm{ad}})_{\underline{c}_{e}})$ is of rank $1$ as a $\mathbb{Z}[v, v^{-1}]$-module. This means we have that the classes of the images of $\varepsilon_{n,\mathcal{L}} * \mathbb{K}_{\mathcal{L}}$ and $\varepsilon_{n,\mathcal{L}}$ under the cell quotient map to $D^b_{\mathfrak{C},\mathcal{L}}(G_{\mathrm{ad}})_{\underline{c}_{e}}$ are scalar multiples, so $[\varepsilon_{n,\mathcal{L}} * \mathbb{K}_{\mathcal{L}}] - q'(v)[\varepsilon_{n,\mathcal{L}}]$ for some $q'(v) \in \mathbb{Z}[v, v^{-1}]$.
    
    By \cite{BT}, $[\mathbb{K}_{\mathcal{L}}] = (v^2 - 1)^{\mathrm{rank}(T)}[\hat{\mathcal{T}}_{\mathcal{L}}]$ in the full Grothendieck group $K_0(\hat{\mathcal{P}}_{\mathcal{L}})$. Note that in the corresponding top cell subquotient module for the monodromic Hecke algebra $K_0(\mathcal{H}_{\mathcal{L}}^{(1)})$, the equation $[\hat{\mathcal{T}}_{\mathcal{L}}] - q_{\mathcal{L}}(v)[\hat{\delta}_{\mathcal{L}}]$ holds by \Cref{lem:yuntilting}. This means that $q'(v) = (v^2 - 1)^{\mathrm{rank}(T)}q_{\mathcal{L}}(v)$ is the only value for which $[\varepsilon_{n,\mathcal{L}} * {\mathbb{K}_{\mathcal{L}}}] - q'(v)[\varepsilon_{n,\mathcal{L}}]$ lies in a lower cell submodule of $K_0(D^b_{\mathfrak{C},\mathcal{L}}(G_{\mathrm{ad}}))$, and therefore $[\varepsilon_{n,\mathcal{L}} * \mathbb{K}_{\mathcal{L}}] = (v^2 - 1)^{\mathrm{rank}(T)}q_{\mathcal{L}}(v)[\varepsilon_{n,\mathcal{L}}]$.
\end{proof}

\subsubsection{Convolution with $\mathbb{K}_{\mathcal{L}}$ for Kazhdan-Laumon objects}

\begin{lemma}\label{lem:csiso}
    For any $s \in S$, the map
    \begin{align}
        c_s * \mathbb{K}_{\mathcal{L}}: {}_{\mathcal{L}}\hat{\nabla}(s)_{s\mathcal{L}} * {}_{s\mathcal{L}}\hat{\nabla}(s)_{\mathcal{L}} * \mathbb{K}_{\mathcal{L}} \to \mathbb{K}_{\mathcal{L}}
    \end{align}
    is an isomorphism.
\end{lemma}

\begin{proof}
    It is enough to show that the corresponding map $\tilde{c}_s : {}_{s\mathcal{L}}\hat{\nabla}(s)_{\mathcal{L}} * \mathbb{K}_{\mathcal{L}} \to {}_{s\mathcal{L}}\hat{\Delta}(s)_{\mathcal{L}} * \mathbb{K}_{\mathcal{L}}$ (obtained by convolving $c_s$ with ${}_{s\mathcal{L}}\hat{\Delta}(s)_{\mathcal{L}}$) is an isomorphism. If $s\mathcal{L} \neq \mathcal{L}$, then by Lemma 3.6 in \cite{lusztig_yun_2020}, ${}_{s\mathcal{L}}\hat{\nabla}(s)_{\mathcal{L}} \cong {}_{s\mathcal{L}}\hat{\Delta}(s)_{\mathcal{L}}$, and so this becomes immediate. We now consider the case when $s\mathcal{L} = \mathcal{L}$. Note that $\tilde{c}_s = (i' \circ p) * \mathbb{K}_{\mathcal{L}}$ where $i'$ and $p$ are as in the canonical exact sequences
    \[\begin{tikzcd}
        0 \arrow[r] & \hat{\mathrm{IC}}(s)_{\mathcal{L}} \arrow[r, "i"] & {}_{s\mathcal{L}}\hat{\nabla}(s)_{\mathcal{L}} \arrow[r, "p"] & \hat{\mathrm{IC}}(e)_{\mathcal{L}} \arrow[r] & 0,\\
        0 \arrow[r] & \hat{\mathrm{IC}}(e)_{\mathcal{L}} \arrow[r, "i'"] & {}_{s\mathcal{L}}\hat{\Delta}(s)_{\mathcal{L}} \arrow[r, "p'"] & \hat{\mathrm{IC}}(s)_{\mathcal{L}} \arrow[r] & 0.
    \end{tikzcd}\]
    It is then enough to show that $p * \mathbb{K}_{\mathcal{L}}$ and $i' * \mathbb{K}_{\mathcal{L}}$ are each isomorphisms. This follows from the fact that $\hat{\mathrm{IC}}(s)_{\mathcal{L}} * \mathbb{K}_{\mathcal{L}} = 0$. Indeed, $\hat{\mathrm{IC}}(s)_{\mathcal{L}} * \mathbb{K} = 0$ if and only if $\hat{\mathrm{IC}}(s)_{\mathcal{L}} * \hat{\mathcal{T}}_{\mathcal{L}} = 0$, which follows from the fact that its class in the Grothendieck is zero combined with the fact that $\hat{\mathcal{T}}_{\mathcal{L}}$ is convolution-exact, since it is tilting.
\end{proof}

\begin{corollary}\label{cor:thetacomp}
    For any $A \in \mathcal{A}_{w,\mathbb{F}_q}$, the object $A * \mathbb{K}_{\mathcal{L}}$ of $\mathcal{A}_{w,\mathbb{F}_q}$ has the property that for any $y, z \in W$, the composition
    \begin{align}
        \theta_{z^{-1},zy} \circ (\Phi_{z^{-1}}^\circ\theta_{z,y}) : A_y * \mathbb{K}_{\mathcal{L}} \to A_y * \mathbb{K}_{\mathcal{L}}
    \end{align}
    is an isomorphism.
\end{corollary}

\begin{proof}
    By \Cref{prop:ftconv}, the morphism in question can be written, up to Tate twist, as a morphism from $A_y * {}_{y\mathcal{L}}\hat{\Delta}(y)_{\mathcal{L}} * \mathbb{K}_{\mathcal{L}} * {}_{\mathcal{L}}\hat{\nabla}(y^{-1})_{y\mathcal{L}} * {}_{y\mathcal{L}}\hat{\nabla}(z^{-1})_{zy\mathcal{L}} * {}_{zy\mathcal{L}}\hat{\nabla}(z)_{y\mathcal{L}}$ to $A_y *{}_{y\mathcal{L}}\hat{\Delta}(y)_{\mathcal{L}} * \mathbb{K}_{\mathcal{L}} * {}_{\mathcal{L}}\hat{\nabla}(y^{-1})_{y\mathcal{L}}$. In particular, it is a Tate twist of the morphism

    \[\begin{tikzcd}
        A_y * {}_{y\mathcal{L}}\hat{\Delta}(y)_{\mathcal{L}} * \mathbb{K}_{\mathcal{L}} * {}_{\mathcal{L}}\hat{\nabla}(y^{-1})_{y\mathcal{L}} * {}_{y\mathcal{L}}\hat{\nabla}(z^{-1})_{zy\mathcal{L}} * {}_{zy\mathcal{L}}\hat{\nabla}(z)_{y\mathcal{L}}\arrow[d]\\
        A_y * {}_{y\mathcal{L}}\hat{\nabla}(z^{-1})_{zy\mathcal{L}} * {}_{zy\mathcal{L}}\hat{\Delta}(zy)_{\mathcal{L}} * \mathbb{K}_{\mathcal{L}} * {}_{zy\mathcal{L}}\hat{\nabla}(y^{-1}z^{-1})_{zy\mathcal{L}} * {}_{zy\mathcal{L}}\hat{\nabla}(z)_{y\mathcal{L}}\arrow[d, "\theta_{z,y} * {}_{zy\mathcal{L}}\hat{\Delta}(zy)_{\mathcal{L}} * \mathbb{K}_{\mathcal{L}} * {}_{\mathcal{L}}\hat{\nabla}(y^{-1}z^{-1})_{zy\mathcal{L}} * {}_{zy\mathcal{L}}\hat{\nabla}(z)_{y\mathcal{L}}"]\\
        A_{zy} * {}_{zy\mathcal{L}}\hat{\Delta}(zy)_{\mathcal{L}} * \mathbb{K}_{\mathcal{L}} * {}_{\mathcal{L}}\hat{\nabla}(y^{-1}z^{-1})_{zy\mathcal{L}} * {}_{zy\mathcal{L}}\hat{\nabla}(z)_{y\mathcal{L}}\arrow[d]\\
        A_{zy} * {}_{zy\mathcal{L}} \hat{\nabla}(z)_{y\mathcal{L}} * {}_{y\mathcal{L}}\hat{\Delta}(y)_{\mathcal{L}} * \mathbb{K}_{\mathcal{L}} * {}_{\mathcal{L}}\hat{\nabla}(y^{-1})_{y\mathcal{L}}\arrow[d, "\theta_{z^{-1},zy} * {}_{y\mathcal{L}}\hat{\Delta}(y)_{\mathcal{L}} * \mathbb{K}_{\mathcal{L}} * {}_{\mathcal{L}}\hat{\nabla}(y^{-1})_{y\mathcal{L}}"]\\
        A_y * {}_{y\mathcal{L}}\hat{\Delta}(y)_{\mathcal{L}} * \mathbb{K}_{\mathcal{L}} * {}_{\mathcal{L}}\hat{\nabla}(y^{-1})_{y\mathcal{L}},
    \end{tikzcd}\]
    where the unlabeled arrows are the isomorphisms given by the central structure on $\mathbb{K}_{\mathcal{L}}$; these extend to similar central morphisms for conjugates of of $\mathbb{K}_{\mathcal{L}}$ by standard/costandard sheaves by the same argument as in Lemma 11.12 of \cite{lusztig_yun_2020}.

    We note that since the second and third morphisms above clearly commute with these central morphisms, the above morphism agrees with the composition
    \[\begin{tikzcd}
        A_y * {}_{y\mathcal{L}}\hat{\Delta}(y)_{\mathcal{L}} * \mathbb{K}_{\mathcal{L}} * {}_{\mathcal{L}}\hat{\nabla}(y^{-1})_{y\mathcal{L}} * {}_{y\mathcal{L}} \hat{\nabla}(z^{-1})_{zy\mathcal{L}} * {}_{zy\mathcal{L}} \hat{\nabla}(z)_{y\mathcal{L}} \arrow[d]\\
        A_y * {}_{y\mathcal{L}}\hat{\nabla}(z^{-1})_{zy\mathcal{L}} * {}_{zy\mathcal{L}}\hat{\nabla}(z)_{y\mathcal{L}} * {}_{y\mathcal{L}}\hat{\Delta}(y)_{\mathcal{L}} * \mathbb{K}_{\mathcal{L}} * {}_{\mathcal{L}}\hat{\nabla}(y^{-1})_{y\mathcal{L}} \arrow[d, "\theta_{z,y} * {}_{zy\mathcal{L}}\hat{\nabla}(z)_{y\mathcal{L}} * {}_{y\mathcal{L}}\hat{\Delta}(y)_{\mathcal{L}} * \mathbb{K}_{\mathcal{L}} * {}_{\mathcal{L}}\hat{\nabla}(y^{-1})_{y\mathcal{L}}"]\\
        A_y * {}_{zy\mathcal{L}}\hat{\nabla}(z)_{y\mathcal{L}} * {}_{y\mathcal{L}}\hat{\Delta}(y)_{\mathcal{L}} * \mathbb{K}_{\mathcal{L}} * {}_{\mathcal{L}}\hat{\nabla}(y^{-1})_{y\mathcal{L}} \arrow[d, "\theta_{z^{-1},zy} * {}_{y\mathcal{L}}\hat{\Delta}(y)_{\mathcal{L}} * \mathbb{K}_{\mathcal{L}} * {}_{\mathcal{L}}\hat{\nabla}(y^{-1})_{y\mathcal{L}}"]\\
        A_y * {}_{y\mathcal{L}}\hat{\Delta}(y)_{\mathcal{L}} * \mathbb{K}_{\mathcal{L}} * {}_{\mathcal{L}}\hat{\nabla}(y^{-1})_{y\mathcal{L}},
    \end{tikzcd}
    \]
   where the first morphism is again the isomorphism coming from centrality of $\mathbb{K}_{\mathcal{L}}$ The composition of the last two morphisms in the sequence above must, by the definition of Kazhdan-Laumon categories, be equal to the morphism 
   \[A_y * c_z * {}_{y\mathcal{L}}\hat{\Delta}(y)_{\mathcal{L}} * \mathbb{K}_{\mathcal{L}} * {}_{\mathcal{L}}\hat{\nabla}(y^{-1})_{y\mathcal{L}},\]
   where $c_z$ is the morphism ${}_{y\mathcal{L}}\hat{\nabla}(z^{-1})_{zy\mathcal{L}} * {}_{zy\mathcal{L}}\hat{\nabla}(z)_{y\mathcal{L}} \to \mathrm{Id}$ which is obtained by applying the morphisms $c_s$ successively for every simple reflection $s$ in a reduced expression for $z$. By the functoriality of the central morphisms discussed above, they also commute with this $c_z$, and so the entire composition above actually agrees with the morphism
   \[A_y * {}_{y\mathcal{L}}\hat{\Delta}(y)_{\mathcal{L}} * \mathbb{K}_{\mathcal{L}} * {}_{\mathcal{L}}\hat{\nabla}(y^{-1})_{y\mathcal{L}} * c_z,\]
 
   By our inductive definition of $c_z$ along with the same argument as in Lemma \ref{lem:csiso}, $\mathbb{K}_{\mathcal{L}} * {}_{\mathcal{L}}\hat{\nabla}(y^{-1})_{y\mathcal{L}} * c_z$ is an isomorphism, and therefore $\theta_{z^{-1},zy} \circ (\Phi_{z^{-1}}^\circ \theta_{z,y})$ must be, too.
\end{proof}

\begin{prop}\label{prop:akfin}
    For any $A \in \mathcal{A}_{w, \mathbb{F}_q}$, $A * \mathbb{K}_{\mathcal{L}}$ has finite projective dimension.
\end{prop}

\begin{proof}
    We can forget the $w$-twisted Weil structure on $A * \mathbb{K}_{\mathcal{L}}$ and consider it as an object in $\mathcal{A}$. In $\mathcal{A}$, the adjunction $(j_{e!}, j_{e}^*)$ gives a morphism
    \begin{align}
        a : j_{e!}j_e^*(A * \mathbb{K}_{\mathcal{L}}) \to A * \mathbb{K}_{\mathcal{L}}.
    \end{align}
    We claim that this is an isomorphism.

    It is enough to show that the component morphisms $a_y : j_y^*j_{e!}j_e^*(A * \mathbb{K}_{\mathcal{L}}) \to j_y^*(A * \mathbb{K}_{\mathcal{L}})$ are each isomorphisms. By definition, these are the structure morphisms $\theta_{y,e} : \Phi_y^\circ(A * \mathbb{K}_{\mathcal{L}})_e \to (A * \mathbb{K}_{\mathcal{L}})_y$.

    By \Cref{cor:thetacomp}, the morphisms 
    \begin{align}
        \theta_{y^{-1},y}\circ (\Phi_{y^{-1}}^\circ \theta_{y,e}) & : \Phi_{y{-1}}^\circ \Phi_y^\circ A_e \to A_e\\
        \theta_{y,e} \circ (\Phi_{y^{-1}}^\circ\theta_{y^{-1},y}) & : \Phi_{y}^\circ\Phi_{y^{-1}}^\circ A_{y} \to A_y
    \end{align}
    are both isomorphisms. This tells us that $\theta_{y,e}$ is both a monomorphism and an epimorphism.
    
    This means $a : j_{e!}j_e^*(A * \mathbb{K}_{\mathcal{L}}) \to A * \mathbb{K}_{\mathcal{L}}$ is an isomorphism as we claimed. The proposition then follows since all objects in the image of $j_{e!}$ have finite cohomological dimension by \Cref{prop:adjunction}.
\end{proof}

\subsection{Completing the proof of Theorem \ref{thm:mainthm}}
\subsubsection{Proof of Theorem \ref{thm:mainthm} by the action of the full twist on cells}

\begin{lemma}\label{lem:ftpolyzero}
    For any $a \in K_0(\mathcal{A}_{w, \mathbb{F}_q}^{\mathcal{L}})$, there exists some $r \geq 1$ for which
    \begin{align}
        P^r(\mathrm{FT}_{\mathcal{L}}, v) \cdot a = 0
    \end{align}
    in $K_0(\mathcal{A}_{w, \mathbb{F}_q}^{\mathcal{L}})$, where $P^r$ is the polynomial for which $P^r(x, v) = P(x, v)^r$ for any $x$.

    Further, if $a \in K_0(\mathcal{A}_{w, \mathbb{F}_q}^{\mathcal{L}})_{< \underline{c}_{e}}$, then
    \begin{align}
        \tilde{P}^r(\mathrm{FT}_{\mathcal{L}}, v) \cdot a = 0.
    \end{align}
\end{lemma}

\begin{proof}
    By and Propositions \ref{prop:dgequiv} and \ref{prop:dgacts}, the category $D^b_{\mathfrak{C},\mathcal{L}}(G_{\mathrm{ad}})$ acts on $D^b_{\mathrm{dg}}(\mathcal{A}_{w,\mathbb{F}_q}^{\mathcal{L}})$ in a way which respects distinguished triangles, therefore giving an action of the $\mathbb{Z}[v, v^{-1}]$-algebra $D^b_{\mathfrak{C},\mathcal{L}}(G_{\mathrm{ad}})$ on $K_0(\mathcal{A}_{w,\mathbb{F}_q})$. It is then enough to show that there exists $r$ for which $P^r([\mathrm{FT}_{\mathcal{L}}], v) = 0$ in $K_0(D^b_{\mathfrak{C},\mathcal{L}}(G_{\mathrm{ad}}))$ and that $\tilde{P}^r([\mathrm{FT}_{\mathcal{L}}], v) \cdot b = 0$ for any $b \in K_0(D^b_{\mathfrak{C},\mathcal{L}}(G_{\mathrm{ad}})_{< \underline{c}_e})$.

    Indeed, since by Proposition \ref{prop:celleigvals} and Proposition \ref{prop:cellsum}, the eigenvalues of the map $[\mathrm{FT}_{\mathcal{L}}] * -$ on $K_0(D^b_{\mathfrak{C},\mathcal{L}}(G_{\mathrm{ad}}))$ are each of the form $v^{2i}$ for $0 \leq i \leq \ell(w_0)$ and of the form $v^{2i}$ for $1 \leq i \leq \ell(w_0)$ on $K_0(D^b_{\mathfrak{C},\mathcal{L}}(G_{\mathrm{ad}})_{< \underline{c}_{e}})$, we must only choose $r$ to be the maximum multiplicity occuring in the characteristic polynomial of $[\mathrm{FT}_{\mathcal{L}}] * -$ on $K_0(D^b_{\mathfrak{C},\mathcal{L}}(G_{\mathrm{ad}}))$, since each degree $1$ term of this characteristic polynomial is a factor of $P(x, v)$ (resp. $\tilde{P}(x, v)$) by definition of the latter. Choosing $r$ in this way, we get that the polynomials in the lemma vanish, as desired.
\end{proof}

The following two lemmas will be used in the proofs of Corollary \ref{cor:lowcellpoly} and Theorem \ref{thm:polyconj}.

\begin{lemma}\label{lem:euclidean}
    Suppose $A$ is a $\mathbb{Z}[v]$-module equipped with a $\mathbb{Z}[v]$-linear endomorphism $T : A \to A$. Let $Q(x, v) \in \mathbb{Z}[x, v]$ be any polynomial. Then if $a \in A$ satisfies
    \begin{align}
        Q(T, v) \cdot a & = 0,
    \end{align}
    then $Q(1, v) a$ lies in the $\mathbb{Z}[T, v]$-span of $(T - 1)\cdot a$.
\end{lemma}
\begin{proof}
    The polynomial $Q(x, v) - Q(1, v) \in \mathbb{Z}[x, v]$ lies in the ideal $(x - 1)$, so we can write
    \begin{align}
        Q(1, v) & = Q(x, v) - \tilde{Q}(x, v)(x - 1)
    \end{align}
    for some polynomial $\tilde{Q}(x, v) \in \mathbb{Z}[x, v]$. Now setting $x = T$ and applying both sides to $a$, we get
    \begin{align}
        Q(1, v)a & = Q(T, v)\cdot a - \tilde{Q}(T, v)(T - 1)\cdot a\\
        & = -\tilde{Q}(T, v)(T - 1)\cdot a,
    \end{align}
    which is in the $\mathbb{Z}[T, v]$-span of $(T - 1)\cdot a$.
\end{proof}

\begin{lemma}\label{lem:ptildezero}
    If $a \in K_0(\mathcal{A}_{w, \mathbb{F}_q}^\mathcal{L})$ is such that $\tilde{P}(\iota^2, v)^k \cdot a = 0$ for any $k \in \mathbb{Z}_{\geq 0}$, then $a \in \vfploc$.
\end{lemma}

\begin{proof}
    Suppose  $a$ satisfies $\tilde{P}(\iota^2, v)^k \cdot a = 0$. By \Cref{thm:i2min1}, we also have $(\iota^2 - 1)^{n} \cdot a \in \vfp$. We claim that this implies that
    \[(\iota^2 - 1)^{n-j} \cdot a \in \vfploc\] for all $0 \leq j \leq n$. Indeed, suppose for induction that \[(\iota^2 - 1)^{n-j+1}\cdot a \in \vfploc,\] and let $a' = (\iota^2 - 1)^{n-j} \cdot a$. Then \[(\iota^2 - 1)\cdot a' \in \vfploc\] and $\tilde{P}(\iota^2, v)^k \cdot a' = 0$, so by Lemma \ref{lem:euclidean} (applied to $T = \iota^2$, $Q(x, v) = \tilde{P}(x, v)^k$), we have that $p(v)^k\cdot a' = \tilde{P}(1, v)^k\cdot a'$ is in the $\mathbb{Z}[\iota^2, v, v^{-1}]$-span of $(\iota^2 - 1) \cdot a'$, and therefore lies in $\vfploc$. Dividing by $p(v)^k$ gives $a' \in \vfploc$, and so we can proceed by induction until $j = n$ where we conclude that $a \in \vfploc$.
\end{proof}

\begin{corollary}\label{cor:lowcellpoly}
For any $a \in K_0(D^b(\mathcal{A}_{w, \mathbb{F}_q}^{\mathcal{L}}))_{< \underline{c}_e}$, $a \in \vfploc$.
\end{corollary}

\begin{proof}
    Let $a \in K_0(\mathcal{A}_{w,\mathbb{F}_q}^{\mathcal{L}})_{< \underline{c}_{e}}$. By Lemmas \ref{lem:i2ft} and \ref{lem:ftpolyzero}, $\tilde{P}(\iota^2, v)^r \cdot a = 0$. Applying Lemma \ref{lem:ptildezero}, we conclude that $a \in \vfploc$.
\end{proof}

We now conclude the proof of \Cref{thm:mainthm}, which states that for any character sheaf $\mathcal{L}$ of $T$ and element $w \in W$, the localization of the $\mathbb{Z}[v, v^{-1}]$-module $K_0(\AwL)$ at $p(v)$ is spanned by classes of objects of finite projective dimension in $\AwL$.

\begin{proof}[Proof of \Cref{thm:mainthm}]
    Let $A \in \mathcal{A}_{w, \mathbb{F}_q}^{\mathcal{L}}$. By Proposition \ref{prop:tiltingid},
    \begin{align}\label{eqn:diff}
        [A * \mathbb{K}_{\mathcal{L}}] - (v^2 - 1)^{\mathrm{rank}(T)}q_{\mathcal{L}}(v)[A]
    \end{align}
    is an element of $K_0(D^b(\mathcal{A}_{w, \mathbb{F}_q}^{\mathcal{L}}))_{< \underline{c}_e}$ 

    By Corollary \ref{cor:lowcellpoly}, this means $(v^2 - 1)^{\mathrm{rank}(T)}([A * \mathbb{K}_{\mathcal{L}}] - q_{\mathcal{L}}(v)[A]) \in \vfploc$.

    By Proposition \ref{prop:akfin}, $A * \mathbb{K}_{\mathcal{L}}$ itself has finite projective dimension, and so in combination with equation (\ref{eqn:diff}), this means $(v^2 - 1)^{\mathrm{rank}(T)}q_{\mathcal{L}}(v)p(v)[A] \in \vfploc$. Note that by \Cref{def:qv} each degree $1$ factor of $q_{\mathcal{L}}(v)$ is also a factor of $p(v)$, so we can divide by $q_{\mathcal{L}}(v)$ in the localization $\vfploc$ to get $[A] \in \vfploc$, completing the proof.
\end{proof}

\section{Polishchuk's rationality conjecture}\label{sec:final}

\subsection{A general study of $K_0(\mathcal{A}_{\mathbb{F}_q})$}
\subsubsection{Polishchuk's description of $K_0(\mathcal{A}_{\mathbb{F}_q})$} 

A crucial tool which we will use in the proof of \Cref{thm:polyconj} is the following description of $K_0(\mathcal{A}_{\mathbb{F}_q})$ provided by Polishchuk.
\begin{theorem}[\cite{P}, Proposition 3.4.1]\label{thm:k0injective}
    The map
    \begin{align}
        K_0(\mathcal{A}_{\mathbb{F}_q}) & \to \bigoplus_{w \in W} K_0(\mathrm{Perv}_{\mathbb{F}_q}(G/U))
    \end{align}
    induced by the functor $\oplus_{w \in W} j_{w}^*$ is injective. Its image is the subset
    \[\{(a_w)_{w \in W} \in K_0(\mathrm{Perv}_{\mathbb{F}_q}(G/U)) ~|~a_{sw} - \Phi_s a_w \in \mathrm{im}(\Phi_s^2 - 1), s \in S, w\in W\}.\]
\end{theorem}
\subsubsection{Recalling \cite{CMFSymplectic}}

In \cite{CMFSymplectic}, we study the subalgebra $\mathrm{KL}(v)$ of endomorphisms of $K_0(G/U)$ generated by the symplectic Fourier transforms $\{\Phi_s\}_{s \in S}$. By \Cref{sec:ft}, this is the same as the subalgebra of $K_0(G/U \times G/U)$ generated under convolution by classes of Kazhdan-Laumon sheaves; we denote the generator of $\mathrm{KL}(v)$ corresponding to $w \in W$ by $\mathsf{a}_w$. In this section, we use the monodromic Hecke algebras $\mathcal{H}_{\mathcal{L}}$ and $\mathcal{H}_{\mathfrak{o}}$ with the standard generators $\tilde{T}_w$ as defined in \cite{lusztig_yun_2020}, see \cite{CMFSymplectic} for a more precise outline of our conventions.

In Section 4 of \cite{CMFSymplectic}, we show that for any character sheaf $\mathcal{L}$ with $W$-orbit $\mathfrak{o}$, there is a surjection $\pi_{\mathcal{L}} : \mathrm{KL}(v) \to \mathcal{H}_{\mathfrak{o}}$. The following result follows from the main result of loc. cit. which explicitly identifies the algebra $\mathrm{KL}(v)$ as a subalgebra of a generic-parameter version of the Yokonuma-Hecke algebra.

\begin{prop}[\cite{CMFSymplectic}]\label{prop:monodrominj}
The following properties are satisfied by the morphisms $\{\pi_{\mathcal{L}}\}_{\mathcal{L}}$.
\begin{enumerate}
    \item If $w_1, w_2 \in W$, then $\pi_{\mathcal{L}}(\mathsf{a}_{w_1}\mathsf{a}_{w_2}) = \pi_{w_2\mathcal{L}}(\mathsf{a}_{w_1})\pi_{\mathcal{L}}(\mathsf{a}_{w_2})$ in $\mathcal{H}_{\mathfrak{o}}$.
    \item If $w \in W_{\mathcal{L}}^\circ$, then $\pi_{\mathcal{L}}(\mathsf{a}_w) = \tilde{T}_{w} \in \mathcal{H}_{\mathfrak{o}}$.
    \item If $s\in S$ is not in $W_{\mathcal{L}}^\circ$, then $\pi_{\mathcal{L}}(\mathsf{a}_s^2) = 1$. 
\end{enumerate}
    Finally, the morphism $\prod_{\mathcal{L}} \pi_{\mathcal{L}}$ is injective, so if $a \in \mathrm{KL}(v)$ is such that $\pi_{\mathcal{L}}(a) = 0$ for all character sheaves $\mathcal{L}$, then $a = 0$. 
\end{prop}

\begin{lemma}\label{lem:w0monodrom}
    For any $\mathcal{L}$ with $W$-orbit $\mathfrak{o}$, the algebra morphism $\pi_{\mathcal{L}} : \mathrm{KL}(v) \to \mathcal{H}_{\mathfrak{o}}$ is such that 
    \begin{align}
        \pi_{\mathcal{L}}(\mathsf{a}_{w_0}^2) & = \tilde{T}_{w_{0,\mathcal{L}}}^2,
    \end{align}
    where $w_{0,\mathcal{L}}$ is the longest element of $W_{\mathcal{L}}^\circ$.
\end{lemma}

\begin{proof}
    Recall that if $s \in S$ is such that $s \not\in W_{\mathcal{L}}^\circ$, then $\pi_{\mathcal{L}}(\mathsf{a}_{s}^2) = 1$ in $\mathcal{H}_{\mathfrak{o}}$. For induction, we claim that if $y \in W$ is such that $\ell(y) + \ell(w_{\mathcal{L},0}) = l(yw_{\mathcal{L},0})$, and if $s \in S$ such that $l(sy) < l(y)$, then $s \not\in W_{y\mathcal{L}}$, i.e. $y^{-1}sy \not\in W_{\mathcal{L}}^\circ$.

    Indeed, if we had $y^{-1}sy \in W_{\mathcal{L}}^\circ$, then we since $w_{0,\mathcal{L}}$ dominates all elements of $W_{\mathcal{L}}^\circ$ in the Bruhat order, we could pick $z \in W$ such that $y^{-1}syz = w_{0,\mathcal{L}}$ with $\ell(y^{-1}sy) + \ell(z) = \ell(w_0,\mathcal{L})$. But then we would have $syz = yw_{0,\mathcal{L}}$ with $\ell(sy) + \ell(z) = \ell(syz) = \ell(yw_{\mathcal{L},0})$. But since $\ell(sy) < \ell(y)$ and $\ell(z) \leq \ell(w_{0,\mathcal{L}})$, this is impossible.

    Now choose some $y$ for which $yw_{0,\mathcal{L}} = w_0$ with $\ell(y) + \ell(w_{0,\mathcal{L}})$. We will show that $\pi_{\mathcal{L}}(\mathsf{a}_{w_0}^2) = \tilde{T}_{w_{0,\mathcal{L}}}^2$ by induction on the length of $y$. Choosing $s \in S$ such that $\ell(sy) < \ell(y)$, this induction hypothesis along with the fact proved in the previous paragraph gives that
    {\allowdisplaybreaks
    \begin{align}
        \pi_{\mathcal{L}}(\mathsf{a}_{w_0}^2) & = \pi_{\mathcal{L}}(\mathsf{a}_{w_{0,\mathcal{L}}}\mathsf{a}_{y^{-1}}\mathsf{a}_y\mathsf{a}_{w_{0,\mathcal{L}}})\\
        & = \pi_{\mathcal{L}}(\mathsf{a}_{w_{0,\mathcal{L}}}\mathsf{a}_{(sy)^{-1}}\mathsf{a}_s^2\mathsf{a}_{sy}\mathsf{a}_{w_{0,\mathcal{L}}})\\
        & = \pi_{sy\mathcal{L}}(\mathsf{a}_{w_{0,\mathcal{L}}}\mathsf{a}_{(sy)^{-1}})\pi_{sy\mathcal{L}}(\mathsf{a}_s^2)\pi_{\mathcal{L}}(\mathsf{a}_{sy}\mathsf{a}_{w_{0,\mathcal{L}}})\\
        & = \pi_{\mathcal{L}}(\mathsf{a}_{w_{0,\mathcal{L}}}\mathsf{a}_{(sy)^{-1}}\mathsf{a}_{sy}\mathsf{a}_{w_{0,\mathcal{L}}})\\
        & = \pi_{\mathcal{L}}(\mathsf{a}_{syw_{0,\mathcal{L}}}^2)\\
        & = \tilde{T}_{w_{0,\mathcal{L}}}^2.
\end{align}}
\end{proof}

\subsubsection{The action of the full twist on $\mathrm{Perv}_{\mathbb{F}_q}(G/U)$}

\begin{prop}\label{prop:i2}
    The endomorphism $\iota^2$ on $K_0(\mathcal{A})$ satisfies $P(\iota^2, v) = 0$. 
    \end{prop}
    
    \begin{proof}
    
    By Theorem \ref{thm:k0injective}, it is enough to show that $P(\Phi_{w_0}^2, v) = 0$ as an endomorphism of the space $K_0(\mathrm{Perv}_{\mathbb{F}_q}(G/U))$.
    
    Now recall that the endomorphism $\Phi_{w_0} : K_0(G/U) \to K_0(G/U)$ agrees with right convolution with the Kazhdan-Laumon sheaf $\overline{K(w_0)}$. Since the convolution $D^b(G/U) \times D^b(G/U \times G/U) \to D^b(G/U)$ is a triangulated functor, it is enough to show that $P([\overline{K(w_0)} * \overline{K(w_0)}], v) = 0$ in $K_0(G/U \times G/U)$.
    
        Letting $1_{\mathcal{L}}$ be the idempotent in $\mathcal{H}_{\mathfrak{o}}$ corresponding to the $\mathcal{L}$-monodromic subalgebra, then in loc. cit. we show that
        \begin{align}
            \pi_{\mathcal{L}}(\mathsf{a}_s) & = \begin{cases}
                -v\tilde{T}_{s}^{-1}1_{\mathcal{L}} & s \in W_{\mathcal{L}}^\circ\\
                -\tilde{T}_s1_{\mathcal{L}} & s \not\in W_{\mathcal{L}}^\circ.
            \end{cases}
        \end{align}
        It is a straightforward calculation from the above to show that $\pi_{\mathcal{L}}([K(w_0) * K(w_0)]) = v^{2\ell(y)}\tilde{T}_{y}^{-2}1_{\mathcal{L}} \in \mathcal{H}_{\mathcal{L}}$, where $y$ is the longest element of $W_{\mathcal{L}}^\circ$. By the main result of \cite{lusztig_yun_2020}, the monodromic Hecke algebra $\mathcal{H}_{\mathcal{L}}$ is isomorphic to $\mathcal{H}_{W_{\mathcal{L}}^\circ}$, with full twists on each side being identified as in \Cref{lem:w0monodrom}. By 5.12.2 of \cite{Lbook} which identifies the eigenvalues of the full twist in the regular representation, we have that $P(-v^{2\ell(y)}\tilde{T}_{y}^2, v) = 0$ in $\mathcal{H}_{\mathcal{L}}$.
    
        Finally, we note that by \Cref{prop:monodrominj}, if a polynomial is satisfied by $\pi_{\mathcal{L}}([\overline{K(w_0)} * \overline{K(w_0)}])$ in each $\mathcal{H}_{\mathfrak{o}}$, then it is also satisfied in $K_0(G/U \times G/U)$; this completes the proof of the proposition.
    \end{proof}

\subsection{Completing the proof of \Cref{thm:polyconj}}

In this section, we complete the proof of \Cref{thm:polyconj}, which states that the localization of the $\mathbb{Z}[v, v^{-1}]$-module $K_0(\AFq)$ at the polynomial
\begin{align}
    p(v) & = \prod_{i=1}^{\ell(w_0)} \left(1 - v^{2i}\right)
\end{align}
is generated by objects of finite projective dimension.

\subsubsection{Proof of \Cref{thm:polyconj}}

Let $a \in K_0(\mathcal{A}_{\mathbb{F}_q})$. We can write $p(v) = \tilde{P}(1, v) = \tilde{P}(x, v) + r(x, v)(x - 1)$ for some $r(x, v) \in \mathbb{Z}[x, v]$. So if
\begin{align}
    a_0 & = \tilde{P}(\iota^2, v)a\\
    a_1 & = r(\iota^2, v)(\iota^2 - 1)a,
\end{align}
then $a_0 + a_1 = p(v)a$, so it suffices to show that $a_0, a_1 \in \vfploc(\mathcal{A}_{\mathbb{F}_q})$.

First we show this for $a_0$. We claim that for any $w \in W$ and $s \in S$, $\Phi_{s}^2((a_0)_w) = (a_0)_w$. Since $a_0 = \tilde{P}(\Phi_{w_0}^2, v)a$, this will follow from the fact that for any $s \in I$, $(\Phi_s^2 - 1)\tilde{P}(\Phi_{w_0}^2, v) = 0$ as an endomorphism of $K_0(G/U)$. By \Cref{prop:monodrominj}, this relation holds if and only if it holds after applying $\pi_{\mathcal{L}}$ for any character sheaf $\mathcal{L}$. By \Cref{lem:w0monodrom}, this reduces to showing that if $\tilde{P}(\tilde{T}_{w_0}^2, v) \neq 0$, then $(\tilde{T}_s^2 - 1)\tilde{P}(\tilde{T}_{w_0}^2, v) = 0$ for any $s \in S$. We can rephrase this as saying that if the full twist acts by the eigenvalue $1$, then so does $\tilde{T}_{s}^2$ for every $s \in S$. Indeed, this follows from the classification of irreducible representations of Hecke algebras, and it is shown directly in 11.5.3 of \cite{P}.

Now note that by \Cref{thm:k0injective}, we have 
\begin{align}\label{eqn:phisquared}
    (a_0)_{sw} - \Phi_s(a_0)_w = (\Phi_s^2 - 1)b
\end{align}
for some $b$. The relation $(\Phi_s + v^2)(\Phi_s^2 - 1) = 0$ from Proposition 6.2.1 of \cite{P} can be rewritten as $\Phi_s(\Phi_s^2 - 1) = -v^2(\Phi_s^2 - 1)$, so when we apply $(\Phi_s^2 - 1)$ to the right-hand side of (\ref{eqn:phisquared}), we get
\begin{align}
    (\Phi_s^2 - 1)^2b & = \Phi_s^2(\Phi_s^2 - 1)b - (\Phi_s^2 - 1)b\\
    & = v^4(\Phi_s^2 - 1)b - (\Phi_s^2 - 1)b\\
    & = (v^4 - 1)(\Phi_s^2 - 1)b.
\end{align}
This means that by applying $(\Phi_s^2 - 1)$ to both sides of (\ref{eqn:phisquared}), and using the fact that $\Phi_s^2((a_0)_{w'}) = (a_0)_{w'}$ for $w' = w$ and $w' = sw$, we then get
    \begin{align}
        0 & = (\Phi_s^2 - 1)((a_0)_{sw} - \Phi_s(a_0)_w)\\
        & = (\Phi_s^2 - 1)^2b\\
        & =(v^4 - 1)(\Phi_{s}^2 - 1)b,\\
        & = (v^4 - 1)((a_0)_{sw} - \Phi_s(a_0)_w),
    \end{align}
so $(v^4 - 1)((a_0)_{sw} - \Phi_s(a_0)_w) = 0$, meaning $(a_0)_{sw} = \Phi_s(a_0)_w$. This means $a_0 = j_{w!}j_w^*(a_0)$ in $K_0(\mathcal{A}_{\mathbb{F}_q})$ for any $w \in W$, and so $a_0 \in \vfploc$.

Now it only remains to observe that $a_1 \in \vfploc$. By Proposition \ref{prop:i2} and the definition of $a_1$, $\tilde{P}(\iota^2, v)a_1 = 0$. By \Cref{lem:ptildezero}, $a_1 \in \vfploc$. Then since $a_0$ and $a_1$ both lie in $\vfploc$ and $p(v)a = a_0 + a_1$, we have that $a \in \vfploc$.

\section{Construction of Kazhdan-Laumon representations}\label{sec:construction}

\subsection{The Grothendieck-Lefschetz pairing}

\subsubsection{The original proposal in \cite{KL}}

In Section 3 of \cite{KL}, the proposed construction of representations is as follows. The authors begin by making Conjecture \ref{conj:kl}, which we now know to be false by Bezrukavnikov and Polishchuk's appendix to \cite{P}. However, for objects of finite projective dimension, one can still define the Grothendieck-Lefschetz-type pairing in the manner they describe

First, they define a Verdier duality functor $\mathbb{D} : \mathcal{A}_\psi \to \mathcal{A}_{\psi^{-1}}$, where $\mathcal{A}_{\psi} = \mathcal{A}$ as we have been using it throughout this paper, while $\mathcal{A}_{\psi^{-1}}$ is the same category but using the additive character $\psi^{-1}$ instead of $\psi$ (where $\psi$ is the additive character we chose in Section \ref{sec:ladic}).

They note that for any $A \in\Aw$ and $B \in (\mathcal{A}_{\psi^{-1}})_{w, \mathbb{F}_q}$ and any $i \in \mathbb{Z}$, the isomorphisms $\psi_A : \mathcal{F}_w\mathrm{Fr}^*A \to A$ and $\psi_B : \mathcal{F}_w\mathrm{Fr}^*B \to B$ give an endomorphism $\psi_{A,B}^i$ of the vector space $\mathrm{Ext}_{\mathcal{A}}^i(A, \mathbb{D}B)$  given by the composition,
\begin{align*}
    & \mathrm{Ext}_{\mathcal{A}}^i(A, \mathbb{D}B) \rightarrow \mathrm{Ext}_{\mathcal{A}}^i(\mathcal{F}_w\mathrm{Fr}^*A, \mathcal{F}_{w}\mathrm{Fr}^*\mathbb{D}B)\\
    & \rightarrow  \mathrm{Ext}_{\mathcal{A}}^i(\mathcal{F}_w A, \mathcal{F}_w \mathbb{D}B) \rightarrow \mathrm{Ext}_{\mathcal{A}}^i(A, \mathbb{D}B) 
\end{align*}
where the first map arises from the morphisms $\psi_A$ and $\psi_B$, the next from the canonical isomorphisms $\mathrm{Fr}^*A \to A$ and $\mathrm{Fr}^*B \to B$, and the last from the fact that $\mathcal{F}_w$ is invertible. This map is also described explicitly in 4.3.1 of \cite{BP}.

We can then define, for $A$ having finite projective dimension and arbitrary $B$, the value
\begin{align}
    \langle [A], [B]\rangle = 
    \sum_{i\in \mathbb{Z}} (-1)^i \mathrm{tr}(\psi^i_{A,B}, \mathrm{Ext}^i_{\mathcal{A}}(A, \mathbb{D}B)).
\end{align}
This is clearly well-defined at the level of Grothendieck groups. We will now explain how to use the result in \Cref{thm:mainthm} to extend Kazhdan and Laumon's pairing, which as of this point is only defined for objects of finite projective dimension, to the full Grothendieck group in the monodromic case.

It is a straightforward computation that for any such $A$ and $B$, we have
\begin{align}
    \langle [A(-\tfrac{1}{2})], [B]\rangle = q^{\frac{1}{2}}\langle [A], [B]\rangle,
\end{align}
so the pairing is $\mathbb{Z}[v, v^{-1}]$-linear where $\mathbb{Z}[v, v^{-1}]$ acts on the target field such that $v$ is multiplication by $q^{\frac{1}{2}}$.

\subsubsection{A pairing on $K_0(\AwL) \otimes \mathbb{C}$}

We can do the same construction on the monodromic Kazhdan-Laumon category $\AwL$ and its Grothendieck group. Then using $\mathbb{Z}[v, v^{-1}]$-linearity, the above definition gives us a pairing which is well-defined on elements of $\vfp$.

Now we note that the polynomial $p(v)$ evaluated at $v = q^{\frac{1}{2}}$ is nonzero, so we can extend this pairing linearly to the localization $\vfploc$. This then gives that the pairing is well-defined on all of $\vfploc \otimes \mathbb{C}$, where in the tensor product we send $v \mapsto q^{\frac{1}{2}}$. But by \Cref{thm:mainthm}, $\vfploc \otimes \mathbb{C}$ is all of $K_0(\AwL) \otimes \mathbb{C}$, so we can indeed define the pairing on this entire vector space; we will now explain how to use this to construct the Kazhdan-Laumon representations.

\subsection{Construction of representations}

As Kazhdan and Laumon explain in \cite{KL}, the category $\mathcal{A}_{w, \mathbb{F}_q}$ is defined so that $K_0(\mathcal{A}_{w, \mathbb{F}_q})$ carries commuting actions of $G(\mathbb{F}_q)$ and $T(w)$, where $T(w)$ is the (usually non-split) torus of $G$ corresponding to $w \in W$, defined by
\begin{align}
    T(w) & = \{t \in T(\overline{\mathbb{F}_q})~|~\mathrm{Fr}^*(t) = w(t)\}.
\end{align}
We note that $K_0(\AwL)$ then also carries commuting actions of $G(\mathbb{F}_q)$ and $T(w)$ where $T(w)$ acts by its character $\theta$ which corresponds to the data of the character sheaf $\mathcal{L}$.

As we explained, the pairing $\langle, \rangle$ is well-defined on $K_0(\AwL)\otimes \mathbb{C}$, and so we can define $K_{w}^{\mathcal{L}}$ to be its kernel. Then the Kazhdan-Laumon representation corresponding to the pairing $(T(w), \theta)$ which was originally sought in \cite{KL} is the vector space $V_{w,\mathcal{L}} = (K_0(\AwL) \otimes \mathbb{C})/K_w^{\mathcal{L}}$. 

In future work, we hope to explicitly decompose this vector space into irreducibles and compute the characters of $V_{w,\mathcal{L}}$ explicitly, generalizing the work which was done in \cite{BP} for quasi-regular characters.

\bibliographystyle{alpha}
\bibliography{bibl}

\begin{thebibliography}{BITV23}

\bibitem[BBP02]{BBP}
Roman Bezrukavnikov, Alexander Braverman, and Leonid Positselskii.
\newblock Gluing of abelian categories and differential operators on the basic affine space.
\newblock {\em J. Inst. Math. Jussieu}, 1(4):543--557, 2002.

\bibitem[BFO12]{BFO}
Roman Bezrukavnikov, Michael Finkelberg, and Victor Ostrik.
\newblock Character {$D$}-modules via {D}rinfeld center of {H}arish-{C}handra bimodules.
\newblock {\em Invent. Math.}, 188(3):589--620, 2012.

\bibitem[BITV23]{BITV}
Roman Bezrukavnikov, Andrei Ionov, Kostiantyn Tolmachov, and Yakov Varshavsky.
\newblock Equivariant derived category of a reductive group as a categorical center, 2023.

\bibitem[BP98]{BP}
Alexander Braverman and Alexander Polishchuk.
\newblock {K}azhdan-{L}aumon representations of finite {C}hevalley groups, character sheaves and some generalization of the {L}efschetz-{V}erdier trace formula, 1998.
\newblock https://arxiv.org/abs/math/9810006.

\bibitem[BT22]{BT}
Roman Bezrukavnikov and Kostiantyn Tolmachov.
\newblock Monodromic model for {K}hovanov-{R}ozansky homology.
\newblock {\em J. Reine Angew. Math.}, 787:79--124, 2022.

\bibitem[BY13]{BY}
Roman Bezrukavnikov and Zhiwei Yun.
\newblock On {K}oszul duality for {K}ac-{M}oody groups.
\newblock {\em Represent. Theory}, 17:1--98, 2013.

\bibitem[BZN15]{BZN}
David Ben-Zvi and David Nadler.
\newblock The character theory of a complex group, 2015.

\bibitem[CYD17]{CYD}
Tsao-Hsien Chen and Alexander Yom~Din.
\newblock A formula for the geometric {J}acquet functor and its character sheaf analogue.
\newblock {\em Geom. Funct. Anal.}, 27(4):772--797, 2017.

\bibitem[Del77]{DCoh}
P.~Deligne.
\newblock {\em Cohomologie \'{e}tale}, volume 569 of {\em Lecture Notes in Math.}
\newblock Springer-Verlag, Berlin, 1977.
\newblock S\'{e}minaire de g\'{e}om\'{e}trie alg\'{e}brique du Bois-Marie SGA $4\frac{1}{2}$.

\bibitem[Gin89]{GinsburgAdmissible}
Victor Ginsburg.
\newblock Admissible modules on a symmetric space.
\newblock {\em Ast\'{e}risque}, 173-174:9--10, 199--255, 1989.
\newblock Orbites unipotentes et repr\'{e}sentations, III.

\bibitem[Gou21]{Gouttard}
Valentin Gouttard.
\newblock {\em {Perverse Monodromic Sheaves}}.
\newblock Theses, {Universit{\'e} Clermont Auvergne}, July 2021.

\bibitem[Gun17]{gunningham}
Sam Gunningham.
\newblock A derived decomposition for equivariant $d$-modules, 2017.

\bibitem[KL88]{KL}
D.~Kazhdan and G.~Laumon.
\newblock Gluing of perverse sheaves and discrete series representation.
\newblock {\em J. Geom. Phys.}, 5(1):63--120, 1988.

\bibitem[Lur17]{ha}
Jacob Lurie.
\newblock Higher algebra.
\newblock Unpublished. Available online at \url{https://people.math.harvard.edu/~lurie/papers/HA.pdf}, 09 2017.

\bibitem[Lus84]{Lbook}
George Lusztig.
\newblock {\em Characters of Reductive Groups over a Finite Field. (AM-107)}.
\newblock Princeton University Press, 1984.

\bibitem[Lus86]{CSIV}
George Lusztig.
\newblock Character sheaves. {IV}.
\newblock {\em Adv. in Math.}, 59(1):1--63, 1986.

\bibitem[LY20]{lusztig_yun_2020}
George Lusztig and Zhiwei Yun.
\newblock Endoscopy for {H}ecke categories, character sheaves and representations.
\newblock {\em Forum of Mathematics, Pi}, 8:e12, 2020.

\bibitem[MF22]{CMFKLCatO}
Calder Morton-Ferguson.
\newblock Kazhdan-{L}aumon {C}ategory $\mathcal{O}$, {B}raverman-{K}azhdan {S}chwartz space, and the semi-infinite flag variety, 2022.

\bibitem[MF23]{CMFSymplectic}
Calder Morton-Ferguson.
\newblock Symplectic {F}ourier-{D}eligne transforms on {G/U} and the algebra of braids and ties, 2023.

\bibitem[Pol01]{P}
Alexander Polishchuk.
\newblock Gluing of perverse sheaves on the basic affine space.
\newblock {\em Selecta Math. (N.S.)}, 7(1):83--147, 2001.
\newblock With an appendix by R. Bezrukavnikov and the author.

\bibitem[Wil03]{W}
Geordie Williamson.
\newblock Mind your {P} and {Q}-symbols: Why the {K}azhdan-{L}usztig basis of the {H}ecke algebra of {T}ype {A} is cellular, 2003.

\bibitem[Yun09]{Z}
Zhiwei Yun.
\newblock Weights of mixed tilting sheaves and geometric {R}ingel duality.
\newblock {\em Selecta Math. (N.S.)}, 14(2):299--320, 2009.

\bibitem[Yun17]{YCh}
Zhiwei Yun.
\newblock Rigidity in the {L}anglands correspondence and applications.
\newblock In {\em Proceedings of the {S}ixth {I}nternational {C}ongress of {C}hinese {M}athematicians. {V}ol. {I}}, volume~36 of {\em Adv. Lect. Math. (ALM)}, pages 199--234. Int. Press, Somerville, MA, 2017.

\end{thebibliography}

\end{document}